\documentclass{amsart}

\usepackage{graphicx}
\usepackage{amssymb}

\newtheorem{defn}{Definition}

\newtheorem{rem}{Remark}
\newtheorem{exa}{Example}
\newtheorem{alg}{Algorithm}

\begin{document}

\title[Table Dirichlet $L$-Series and Prime Zeta Modulo Functions]{Table of Dirichlet $L$-Series and Prime Zeta Modulo Functions for Small Moduli}

\author{Richard J. Mathar}
\urladdr{http://www.mpia-hd.mpg.de/~mathar}
\email{mathar@mpia.de}
\address{Max-Planck Institute of Astronomy, K\"onigstuhl 17, 69117 Heidelberg, Germany}


\date{\today}

\begin{abstract}
The Dirichlet characters of reduced residue systems modulo $m$
are tabulated for moduli
$m \le 195$.
The associated $L$-series
are tabulated for $m\le 14$ and small positive integer argument $s$
accurate to $10^{-50}$, their first derivatives for $m\le 6$.
Restricted summation over primes only defines
Dirichlet Prime $L$-functions which lead to Euler products
(Prime Zeta Modulo functions).
Both are materialized over similar ranges of moduli and arguments.
Formulas and numerical techniques are well known; the aim is to
provide direct access to reference values.
\end{abstract}

\maketitle
\section{Introduction} \label{sec.intro} 
\subsection{Scope} 
Dirichlet $L$-series are a standard modification of the Riemann
series of the $\zeta$-function with the intent to distribute the sum over classes
of a modulo system.
Section \ref{sec.intro} enumerates the group representations; a first use
of this classification is a table of $L$-series and some of their first
derivatives at small integer arguments in Section \ref{sec.L}.

In the same spirit as $L$-series generalize the $\zeta$-series,
multiplying the terms of the Prime Zeta Function by the characters
imprints a periodic texture on these, which can be synthesized
(in the Fourier sense) to distil the Prime Zeta Modulo Functions. Section \ref{sec.P}
recalls numerical techniques and tabulates these for small
moduli and small integer arguments.

If some of the Euler products that arise in growth rate estimators
of number densities are to be factorized over the residue classes
of the primes, taking logarithms proposes to use the
Prime Zeta Modulo Function as a basis. The fundamental numerical examples
are worked out in Section \ref{sec.HL}.

\subsection{Dirichlet Characters} \label{sec.char} 

The Dirichlet characters $\chi_r(n)$ (mod $m$) are shown in
Table
\ref{tab.m2}--\ref{tab.m22}
for small modulus $m$.
The purpose of rolling out such basic information is that it
tags each representation with a unique $r$ for use in all tables further down.

Each table shows one character per line, lines enumerated by representation $r$,
$1\le r\le \varphi(m)$, where $\varphi$ is Euler's totient function.
The principal character $\chi_1$ is the top line.
The residues from 1 to $m$ are indicated in the header row.
The entries are either zero or roots of unity. 
As a shortcut to the notation,
\begin{equation}
u_j \equiv e^{2\pi i j/\varphi(m)}; \quad
\bar u_j \equiv e^{-2\pi i j/\varphi(m)},
\end{equation}
denote a root of unity and its complex conjugation.
The conductor $f$ (smallest induced modulus) is another
column;
it also indicates
whether the character is (im)primitive.
Six-digit sequence numbers to matching periodic
sequences in the Online Encyclopedia of Integer Sequences (OEIS) are included \cite{EIS}.

The cases $m\le 7$ have been tabulated by Apostol \cite{Apostol},
some of the real characters by Davies and Haselgrove \cite{DaviesPRSA264},
all up to $m=10$ by Zucker and McPhedran \cite{ZuckerPRSA464}.

The calculations are based on standard algorithms:
$m$ is decomposed into its prime number factorization, the
character table of each factor is constructed and these representations
are multiplied following the group property (multiplication rule)
for each pair of representations, looping over all
products of representations
\cite{Apostol,SpiraMathComp23,ZuckerJPA9}.

The number of real characters $m>2$ is given
by the terms of the OEIS sequence \cite[A060594]{EIS}.
The number of ``quartic'' chacters for which each element is real or purely imaginary
is \cite[A073103]{EIS}.

\begin{table}
\caption{$\chi_r(n)$, $m=2$, $\varphi(m)=1$.}
\begin{tabular}{r|cc|c|r}
\hline
$r$ & 1 & 2 & $f$\\\hline
1 &1 &0 &1 &A000035\\

\hline
\end{tabular}
\label{tab.m2}
\end{table}

\begin{table}
\caption{$\chi_r(n)$, $m=3$, $\varphi(m)=2$.}
\begin{tabular}{r|ccc|c|r}
\hline
$r$ & 1 & 2 & 3 & $f$\\\hline
1 &1 &1 &0 &1 &A011655\\
2 &1 &-1 &0 &3 &A102283\\

\hline
\end{tabular}
\label{tab.m3}
\end{table}

\begin{table}
\caption{$\chi_r(n)$, $m=4$, $\varphi(m)=2$.}
\begin{tabular}{r|cccc|c|r}
\hline
$r$ & 1 & 2 & 3 & 4 & $f$\\\hline
1 &1 &0 &1 &0 &1 &A000035\\
2 &1 &0 &-1 &0 &4 &A101455\\

\hline
\end{tabular}
\label{tab.m4}
\end{table}

\begin{table}
\caption{$\chi_r(n)$, $m=5$, $\varphi(m)=4$.}
\begin{tabular}{r|ccccc|c|r}
\hline
$r$ & 1 & 2 & 3 & 4 & 5 & $f$\\\hline
1 &1 &1 &1 &1 &0 &1 &A011558\\
2 &1 &$i$ &$-i$ &-1 &0 &5 &\\
3 &1 &-1 &-1 &1 &0 &5 &A080891\\
4 &1 &$-i$ &$i$ &-1 &0 &5 &\\

\hline
\end{tabular}
\label{tab.m5}
\end{table}
\clearpage

\begin{table}
\caption{$\chi_r(n)$, $m=6$, $\varphi(m)=2$.}
\begin{tabular}{r|cccccc|c|r}
\hline
$r$ & 1 & 2 & 3 & 4 & 5 & 6 & $f$\\\hline
1 &1 &0 &0 &0 &1 &0 &1 &A120325\\
2 &1 &0 &0 &0 &-1 &0 &3 &A134667\\

\hline
\end{tabular}
\label{tab.m6}
\end{table}

\begin{table}
\caption{$\chi_r(n)$, $m=7$, $\varphi(m)=6$.}
\begin{tabular}{r|ccccccc|c|r}
\hline
$r$ & 1 & 2 & 3 & 4 & 5 & 6 & 7 & $f$\\\hline
1 &1 &1 &1 &1 &1 &1 &0 &1 &A109720\\
2 &1 &$u_{2}$ &$u_{1}$ &$\bar u_{2}$ &$\bar u_{1}$ &-1 &0 &7 &\\
3 &1 &$\bar u_{2}$ &$u_{2}$ &$u_{2}$ &$\bar u_{2}$ &1 &0 &7 &\\
4 &1 &1 &-1 &1 &-1 &-1 &0 &7 &A175629\\
5 &1 &$u_{2}$ &$\bar u_{2}$ &$\bar u_{2}$ &$u_{2}$ &1 &0 &7 &\\
6 &1 &$\bar u_{2}$ &$\bar u_{1}$ &$u_{2}$ &$u_{1}$ &-1 &0 &7 &\\

\hline
\end{tabular}
\label{tab.m7}
\end{table}

\begin{table}
\caption{$\chi_r(n)$, $m=8$, $\varphi(m)=4$.}
\begin{tabular}{r|cccccccc|c|r}
\hline
$r$ & 1 & 2 & 3 & 4 & 5 & 6 & 7 & 8 & $f$\\\hline
1 &1 &0 &1 &0 &1 &0 &1 &0 &1 &A000035\\
2 &1 &0 &-1 &0 &-1 &0 &1 &0 &8 &A091337\\
3 &1 &0 &-1 &0 &1 &0 &-1 &0 &4 &A101455\\
4 &1 &0 &1 &0 &-1 &0 &-1 &0 &8 &\\

\hline
\end{tabular}
\label{tab.m8}
\end{table}

\begin{table}
\caption{$\chi_r(n)$, $m=9$, $\varphi(m)=6$.}
\begin{tabular}{r|ccccccccc|c|r}
\hline
$r$ & 1 & 2 & 3 & 4 & 5 & 6 & 7 & 8 & 9 & $f$\\\hline
1 &1 &1 &0 &1 &1 &0 &1 &1 &0 &1 &A011655\\
2 &1 &$u_{1}$ &0 &$u_{2}$ &$\bar u_{1}$ &0 &$\bar u_{2}$ &-1 &0 &9 &\\
3 &1 &$u_{2}$ &0 &$\bar u_{2}$ &$\bar u_{2}$ &0 &$u_{2}$ &1 &0 &9 &\\
4 &1 &-1 &0 &1 &-1 &0 &1 &-1 &0 &3 &A102283\\
5 &1 &$\bar u_{2}$ &0 &$u_{2}$ &$u_{2}$ &0 &$\bar u_{2}$ &1 &0 &9 &\\
6 &1 &$\bar u_{1}$ &0 &$\bar u_{2}$ &$u_{1}$ &0 &$u_{2}$ &-1 &0 &9 &\\

\hline
\end{tabular}
\label{tab.m9}
\end{table}

\begin{table}
\caption{$\chi_r(n)$, $m=10$, $\varphi(m)=4$.}
\begin{tabular}{r|cccccccccc|c|r}
\hline
$r$ & 1 & 2 & 3 & 4 & 5 & 6 & 7 & 8 & 9 & 10 & $f$\\\hline
1 &1 &0 &1 &0 &0 &0 &1 &0 &1 &0 &1 &\\
2 &1 &0 &$-i$ &0 &0 &0 &$i$ &0 &-1 &0 &5 &\\
3 &1 &0 &-1 &0 &0 &0 &-1 &0 &1 &0 &5 &\\
4 &1 &0 &$i$ &0 &0 &0 &$-i$ &0 &-1 &0 &5 &\\

\hline
\end{tabular}
\label{tab.m10}
\end{table}

\begin{table}
\caption{$\chi_r(n)$, $m=11$, $\varphi(m)=10$.}
\begin{tabular}{r|ccccccccccc|c|r}
\hline
$r$ & 1 & 2 & 3 & 4 & 5 & 6 & 7 & 8 & 9 & 10 & 11 & $f$\\\hline
1 &1 &1 &1 &1 &1 &1 &1 &1 &1 &1 &0 &1 &A145568\\
2 &1 &$u_{1}$ &$\bar u_{2}$ &$u_{2}$ &$u_{4}$ &$\bar u_{1}$ &$\bar u_{3}$ &$u_{3}$ &$\bar u_{4}$ &-1 &0 &11 &\\
3 &1 &$u_{2}$ &$\bar u_{4}$ &$u_{4}$ &$\bar u_{2}$ &$\bar u_{2}$ &$u_{4}$ &$\bar u_{4}$ &$u_{2}$ &1 &0 &11 &\\
4 &1 &$u_{3}$ &$u_{4}$ &$\bar u_{4}$ &$u_{2}$ &$\bar u_{3}$ &$u_{1}$ &$\bar u_{1}$ &$\bar u_{2}$ &-1 &0 &11 &\\
5 &1 &$u_{4}$ &$u_{2}$ &$\bar u_{2}$ &$\bar u_{4}$ &$\bar u_{4}$ &$\bar u_{2}$ &$u_{2}$ &$u_{4}$ &1 &0 &11 &\\
6 &1 &-1 &1 &1 &1 &-1 &-1 &-1 &1 &-1 &0 &11 &A011582\\
7 &1 &$\bar u_{4}$ &$\bar u_{2}$ &$u_{2}$ &$u_{4}$ &$u_{4}$ &$u_{2}$ &$\bar u_{2}$ &$\bar u_{4}$ &1 &0 &11 &\\
8 &1 &$\bar u_{3}$ &$\bar u_{4}$ &$u_{4}$ &$\bar u_{2}$ &$u_{3}$ &$\bar u_{1}$ &$u_{1}$ &$u_{2}$ &-1 &0 &11 &\\
9 &1 &$\bar u_{2}$ &$u_{4}$ &$\bar u_{4}$ &$u_{2}$ &$u_{2}$ &$\bar u_{4}$ &$u_{4}$ &$\bar u_{2}$ &1 &0 &11 &\\
10 &1 &$\bar u_{1}$ &$u_{2}$ &$\bar u_{2}$ &$\bar u_{4}$ &$u_{1}$ &$u_{3}$ &$\bar u_{3}$ &$u_{4}$ &-1 &0 &11 &\\

\hline
\end{tabular}
\label{tab.m11}
\end{table}

\begin{table}
\caption{$\chi_r(n)$, $m=12$, $\varphi(m)=4$.}
\begin{tabular}{r|cccccccccccc|c|r}
\hline
$r$ & 1 & 2 & 3 & 4 & 5 & 6 & 7 & 8 & 9 & 10 & 11 & 12 & $f$\\\hline
1 &1 &0 &0 &0 &1 &0 &1 &0 &0 &0 &1 &0 &1 &A120325\\
2 &1 &0 &0 &0 &-1 &0 &1 &0 &0 &0 &-1 &0 &3 &A134667\\
3 &1 &0 &0 &0 &1 &0 &-1 &0 &0 &0 &-1 &0 &4 &\\
4 &1 &0 &0 &0 &-1 &0 &-1 &0 &0 &0 &1 &0 &12 &A110161\\

\hline
\end{tabular}
\label{tab.m12}
\end{table}

\begin{table}
\caption{$\chi_r(n)$, $m=13$, $\varphi(m)=12$.}
\begin{tabular}{r|ccccccccccccc|c|r}
\hline
$r$ & 1 & 2 & 3 & 4 & 5 & 6 & 7 & 8 & 9 & 10 & 11 & 12 & 13 & $f$\\\hline
1 &1 &1 &1 &1 &1 &1 &1 &1 &1 &1 &1 &1 &0 &1 &\\
2 &1 &$u_{1}$ &$u_{4}$ &$u_{2}$ &$-i$ &$u_{5}$ &$\bar u_{1}$ &$i$ &$\bar u_{4}$ &$\bar u_{2}$ &$\bar u_{5}$ &-1 &0 &13 &\\
3 &1 &$u_{2}$ &$\bar u_{4}$ &$u_{4}$ &-1 &$\bar u_{2}$ &$\bar u_{2}$ &-1 &$u_{4}$ &$\bar u_{4}$ &$u_{2}$ &1 &0 &13 &\\
4 &1 &$i$ &1 &-1 &$i$ &$i$ &$-i$ &$-i$ &1 &-1 &$-i$ &-1 &0 &13 &\\
5 &1 &$u_{4}$ &$u_{4}$ &$\bar u_{4}$ &1 &$\bar u_{4}$ &$\bar u_{4}$ &1 &$\bar u_{4}$ &$u_{4}$ &$u_{4}$ &1 &0 &13 &\\
6 &1 &$u_{5}$ &$\bar u_{4}$ &$\bar u_{2}$ &$-i$ &$u_{1}$ &$\bar u_{5}$ &$i$ &$u_{4}$ &$u_{2}$ &$\bar u_{1}$ &-1 &0 &13 &\\
7 &1 &-1 &1 &1 &-1 &-1 &-1 &-1 &1 &1 &-1 &1 &0 &13 &A011583\\
8 &1 &$\bar u_{5}$ &$u_{4}$ &$u_{2}$ &$i$ &$\bar u_{1}$ &$u_{5}$ &$-i$ &$\bar u_{4}$ &$\bar u_{2}$ &$u_{1}$ &-1 &0 &13 &\\
9 &1 &$\bar u_{4}$ &$\bar u_{4}$ &$u_{4}$ &1 &$u_{4}$ &$u_{4}$ &1 &$u_{4}$ &$\bar u_{4}$ &$\bar u_{4}$ &1 &0 &13 &\\
10 &1 &$-i$ &1 &-1 &$-i$ &$-i$ &$i$ &$i$ &1 &-1 &$i$ &-1 &0 &13 &\\
11 &1 &$\bar u_{2}$ &$u_{4}$ &$\bar u_{4}$ &-1 &$u_{2}$ &$u_{2}$ &-1 &$\bar u_{4}$ &$u_{4}$ &$\bar u_{2}$ &1 &0 &13 &\\
12 &1 &$\bar u_{1}$ &$\bar u_{4}$ &$\bar u_{2}$ &$i$ &$\bar u_{5}$ &$u_{1}$ &$-i$ &$u_{4}$ &$u_{2}$ &$u_{5}$ &-1 &0 &13 &\\

\hline
\end{tabular}
\label{tab.m13}
\end{table}

\begin{table}
\caption{$\chi_r(n)$, $m=14$, $\varphi(m)=6$.}
\begin{tabular}{r|cccccccccccccc|c|r}
\hline
$r$ & 1 & 2 & 3 & 4 & 5 & 6 & 7 & 8 & 9 & 10 & 11 & 12 & 13 & 14 & $f$\\\hline
1 &1 &0 &1 &0 &1 &0 &0 &0 &1 &0 &1 &0 &1 &0 &1 &\\
2 &1 &0 &$u_{1}$ &0 &$\bar u_{1}$ &0 &0 &0 &$u_{2}$ &0 &$\bar u_{2}$ &0 &-1 &0 &7 &\\
3 &1 &0 &$u_{2}$ &0 &$\bar u_{2}$ &0 &0 &0 &$\bar u_{2}$ &0 &$u_{2}$ &0 &1 &0 &7 &\\
4 &1 &0 &-1 &0 &-1 &0 &0 &0 &1 &0 &1 &0 &-1 &0 &7 &\\
5 &1 &0 &$\bar u_{2}$ &0 &$u_{2}$ &0 &0 &0 &$u_{2}$ &0 &$\bar u_{2}$ &0 &1 &0 &7 &\\
6 &1 &0 &$\bar u_{1}$ &0 &$u_{1}$ &0 &0 &0 &$\bar u_{2}$ &0 &$u_{2}$ &0 &-1 &0 &7 &\\

\hline
\end{tabular}
\label{tab.m14}
\end{table}

\begin{table}
\caption{$\chi_r(n)$, $m=15$, $\varphi(m)=8$.}
\begin{tabular}{r|ccccccccccccccc|c|r}
\hline
$r$ & 1 & 2 & 3 & 4 & 5 & 6 & 7 & 8 & 9 & 10 & 11 & 12 & 13 & 14 & 15 & $f$\\\hline
1 &1 &1 &0 &1 &0 &0 &1 &1 &0 &0 &1 &0 &1 &1 &0 &1 &\\
2 &1 &$i$ &0 &-1 &0 &0 &$i$ &$-i$ &0 &0 &1 &0 &$-i$ &-1 &0 &5 &\\
3 &1 &-1 &0 &1 &0 &0 &-1 &-1 &0 &0 &1 &0 &-1 &1 &0 &5 &\\
4 &1 &$-i$ &0 &-1 &0 &0 &$-i$ &$i$ &0 &0 &1 &0 &$i$ &-1 &0 &5 &\\
5 &1 &-1 &0 &1 &0 &0 &1 &-1 &0 &0 &-1 &0 &1 &-1 &0 &3 &\\
6 &1 &$-i$ &0 &-1 &0 &0 &$i$ &$i$ &0 &0 &-1 &0 &$-i$ &1 &0 &15 &\\
7 &1 &1 &0 &1 &0 &0 &-1 &1 &0 &0 &-1 &0 &-1 &-1 &0 &15 &\\
8 &1 &$i$ &0 &-1 &0 &0 &$-i$ &$-i$ &0 &0 &-1 &0 &$i$ &1 &0 &15 &\\

\hline
\end{tabular}
\label{tab.m15}
\end{table}

\begin{table}
\caption{$\chi_r(n)$, $m=16$, $\varphi(m)=8$.}
\begin{tabular}{r|cccccccccccccccc|c|r}
\hline
$r$ & 1 & 2 & 3 & 4 & 5 & 6 & 7 & 8 & 9 & 10 & 11 & 12 & 13 & 14 & 15 & 16 & $f$\\\hline
1 &1 &0 &1 &0 &1 &0 &1 &0 &1 &0 &1 &0 &1 &0 &1 &0 &1 &\\
2 &1 &0 &$i$ &0 &$-i$ &0 &-1 &0 &-1 &0 &$-i$ &0 &$i$ &0 &1 &0 &16 &\\
3 &1 &0 &-1 &0 &-1 &0 &1 &0 &1 &0 &-1 &0 &-1 &0 &1 &0 &8 &\\
4 &1 &0 &$-i$ &0 &$i$ &0 &-1 &0 &-1 &0 &$i$ &0 &$-i$ &0 &1 &0 &16 &\\
5 &1 &0 &-1 &0 &1 &0 &-1 &0 &1 &0 &-1 &0 &1 &0 &-1 &0 &4 &\\
6 &1 &0 &$-i$ &0 &$-i$ &0 &1 &0 &-1 &0 &$i$ &0 &$i$ &0 &-1 &0 &16 &\\
7 &1 &0 &1 &0 &-1 &0 &-1 &0 &1 &0 &1 &0 &-1 &0 &-1 &0 &8 &\\
8 &1 &0 &$i$ &0 &$i$ &0 &1 &0 &-1 &0 &$-i$ &0 &$-i$ &0 &-1 &0 &16 &\\

\hline
\end{tabular}
\label{tab.m16}
\end{table}

\begin{table}
\caption{$\chi_r(n)$, $m=17$, $\varphi(m)=16$.}
\begin{tabular}{r|ccccccccccccccccc|c|r}
\hline
$r$ & 1 & 2 & 3 & 4 & 5 & 6 & 7 & 8 & 9 & 10 & 11 & 12 & 13 & 14 & 15 & 16 & 17 & $f$\\\hline
1 &1 &1 &1 &1 &1 &1 &1 &1 &1 &1 &1 &1 &1 &1 &1 &1 &0 &1 &\\
2 &1 &$\bar u_{2}$ &$u_{1}$ &$-i$ &$u_{5}$ &$\bar u_{1}$ &$\bar u_{5}$ &$\bar u_{6}$ &$u_{2}$ &$u_{3}$ &$u_{7}$ &$\bar u_{3}$ &$i$ &$\bar u_{7}$ &$u_{6}$ &-1 &0 &17 &\\
3 &1 &$-i$ &$u_{2}$ &-1 &$\bar u_{6}$ &$\bar u_{2}$ &$u_{6}$ &$i$ &$i$ &$u_{6}$ &$\bar u_{2}$ &$\bar u_{6}$ &-1 &$u_{2}$ &$-i$ &1 &0 &17 &\\
4 &1 &$\bar u_{6}$ &$u_{3}$ &$i$ &$\bar u_{1}$ &$\bar u_{3}$ &$u_{1}$ &$\bar u_{2}$ &$u_{6}$ &$\bar u_{7}$ &$u_{5}$ &$u_{7}$ &$-i$ &$\bar u_{5}$ &$u_{2}$ &-1 &0 &17 &\\
5 &1 &-1 &$i$ &1 &$i$ &$-i$ &$-i$ &-1 &-1 &$-i$ &$-i$ &$i$ &1 &$i$ &-1 &1 &0 &17 &\\
6 &1 &$u_{6}$ &$u_{5}$ &$-i$ &$\bar u_{7}$ &$\bar u_{5}$ &$u_{7}$ &$u_{2}$ &$\bar u_{6}$ &$\bar u_{1}$ &$u_{3}$ &$u_{1}$ &$i$ &$\bar u_{3}$ &$\bar u_{2}$ &-1 &0 &17 &\\
7 &1 &$i$ &$u_{6}$ &-1 &$\bar u_{2}$ &$\bar u_{6}$ &$u_{2}$ &$-i$ &$-i$ &$u_{2}$ &$\bar u_{6}$ &$\bar u_{2}$ &-1 &$u_{6}$ &$i$ &1 &0 &17 &\\
8 &1 &$u_{2}$ &$u_{7}$ &$i$ &$u_{3}$ &$\bar u_{7}$ &$\bar u_{3}$ &$u_{6}$ &$\bar u_{2}$ &$u_{5}$ &$u_{1}$ &$\bar u_{5}$ &$-i$ &$\bar u_{1}$ &$\bar u_{6}$ &-1 &0 &17 &\\
9 &1 &1 &-1 &1 &-1 &-1 &-1 &1 &1 &-1 &-1 &-1 &1 &-1 &1 &1 &0 &17 &\\
10 &1 &$\bar u_{2}$ &$\bar u_{7}$ &$-i$ &$\bar u_{3}$ &$u_{7}$ &$u_{3}$ &$\bar u_{6}$ &$u_{2}$ &$\bar u_{5}$ &$\bar u_{1}$ &$u_{5}$ &$i$ &$u_{1}$ &$u_{6}$ &-1 &0 &17 &\\
11 &1 &$-i$ &$\bar u_{6}$ &-1 &$u_{2}$ &$u_{6}$ &$\bar u_{2}$ &$i$ &$i$ &$\bar u_{2}$ &$u_{6}$ &$u_{2}$ &-1 &$\bar u_{6}$ &$-i$ &1 &0 &17 &\\
12 &1 &$\bar u_{6}$ &$\bar u_{5}$ &$i$ &$u_{7}$ &$u_{5}$ &$\bar u_{7}$ &$\bar u_{2}$ &$u_{6}$ &$u_{1}$ &$\bar u_{3}$ &$\bar u_{1}$ &$-i$ &$u_{3}$ &$u_{2}$ &-1 &0 &17 &\\
13 &1 &-1 &$-i$ &1 &$-i$ &$i$ &$i$ &-1 &-1 &$i$ &$i$ &$-i$ &1 &$-i$ &-1 &1 &0 &17 &\\
14 &1 &$u_{6}$ &$\bar u_{3}$ &$-i$ &$u_{1}$ &$u_{3}$ &$\bar u_{1}$ &$u_{2}$ &$\bar u_{6}$ &$u_{7}$ &$\bar u_{5}$ &$\bar u_{7}$ &$i$ &$u_{5}$ &$\bar u_{2}$ &-1 &0 &17 &\\
15 &1 &$i$ &$\bar u_{2}$ &-1 &$u_{6}$ &$u_{2}$ &$\bar u_{6}$ &$-i$ &$-i$ &$\bar u_{6}$ &$u_{2}$ &$u_{6}$ &-1 &$\bar u_{2}$ &$i$ &1 &0 &17 &\\
16 &1 &$u_{2}$ &$\bar u_{1}$ &$i$ &$\bar u_{5}$ &$u_{1}$ &$u_{5}$ &$u_{6}$ &$\bar u_{2}$ &$\bar u_{3}$ &$\bar u_{7}$ &$u_{3}$ &$-i$ &$u_{7}$ &$\bar u_{6}$ &-1 &0 &17 &\\

\hline
\end{tabular}
\label{tab.m17}
\end{table}
% \clearpage

\begin{table}
\caption{$\chi_r(n)$, $m=18$, $\varphi(m)=6$.}
\begin{tabular}{r|cccccccccccccccccc|c|r}
\hline
$r$ & 1 & 2 & 3 & 4 & 5 & 6 & 7 & 8 & 9 & 10 & 11 & 12 & 13 & 14 & 15 & 16 & 17 & 18 & $f$\\\hline
1 &1 &0 &0 &0 &1 &0 &1 &0 &0 &0 &1 &0 &1 &0 &0 &0 &1 &0 &1 &\\
2 &1 &0 &0 &0 &$\bar u_{1}$ &0 &$\bar u_{2}$ &0 &0 &0 &$u_{1}$ &0 &$u_{2}$ &0 &0 &0 &-1 &0 &9 &\\
3 &1 &0 &0 &0 &$\bar u_{2}$ &0 &$u_{2}$ &0 &0 &0 &$u_{2}$ &0 &$\bar u_{2}$ &0 &0 &0 &1 &0 &9 &\\
4 &1 &0 &0 &0 &-1 &0 &1 &0 &0 &0 &-1 &0 &1 &0 &0 &0 &-1 &0 &3 &\\
5 &1 &0 &0 &0 &$u_{2}$ &0 &$\bar u_{2}$ &0 &0 &0 &$\bar u_{2}$ &0 &$u_{2}$ &0 &0 &0 &1 &0 &9 &\\
6 &1 &0 &0 &0 &$u_{1}$ &0 &$u_{2}$ &0 &0 &0 &$\bar u_{1}$ &0 &$\bar u_{2}$ &0 &0 &0 &-1 &0 &9 &\\

\hline
\end{tabular}
\label{tab.m18}
\end{table}

\begin{table}
\caption{$\chi_r(n)$, $m=19$, $\varphi(m)=18$.}
\begin{tabular}{r|ccccccccccccccccccc|c|r}
\hline
$r$ & 1 & 2 & 3 & 4 & 5 & 6 & 7 & 8 & 9 & 10 & 11 & 12 & 13 & 14 & 15 & 16 & 17 & 18 & 19 & $f$\\\hline
1 &1 &1 &1 &1 &1 &1 &1 &1 &1 &1 &1 &1 &1 &1 &1 &1 &1 &1 &0 &1 &\\
2 &1 &$u_{1}$ &$\bar u_{5}$ &$u_{2}$ &$\bar u_{2}$ &$\bar u_{4}$ &$u_{6}$ &$u_{3}$ &$u_{8}$ &$\bar u_{1}$ &$\bar u_{6}$ &$\bar u_{3}$ &$u_{5}$ &$u_{7}$ &$\bar u_{7}$ &$u_{4}$ &$\bar u_{8}$ &-1 &0 &19 &\\
3 &1 &$u_{2}$ &$u_{8}$ &$u_{4}$ &$\bar u_{4}$ &$\bar u_{8}$ &$\bar u_{6}$ &$u_{6}$ &$\bar u_{2}$ &$\bar u_{2}$ &$u_{6}$ &$\bar u_{6}$ &$\bar u_{8}$ &$\bar u_{4}$ &$u_{4}$ &$u_{8}$ &$u_{2}$ &1 &0 &19 &\\
4 &1 &$u_{3}$ &$u_{3}$ &$u_{6}$ &$\bar u_{6}$ &$u_{6}$ &1 &-1 &$u_{6}$ &$\bar u_{3}$ &1 &-1 &$\bar u_{3}$ &$u_{3}$ &$\bar u_{3}$ &$\bar u_{6}$ &$\bar u_{6}$ &-1 &0 &19 &\\
5 &1 &$u_{4}$ &$\bar u_{2}$ &$u_{8}$ &$\bar u_{8}$ &$u_{2}$ &$u_{6}$ &$\bar u_{6}$ &$\bar u_{4}$ &$\bar u_{4}$ &$\bar u_{6}$ &$u_{6}$ &$u_{2}$ &$\bar u_{8}$ &$u_{8}$ &$\bar u_{2}$ &$u_{4}$ &1 &0 &19 &\\
6 &1 &$u_{5}$ &$\bar u_{7}$ &$\bar u_{8}$ &$u_{8}$ &$\bar u_{2}$ &$\bar u_{6}$ &$\bar u_{3}$ &$u_{4}$ &$\bar u_{5}$ &$u_{6}$ &$u_{3}$ &$u_{7}$ &$\bar u_{1}$ &$u_{1}$ &$u_{2}$ &$\bar u_{4}$ &-1 &0 &19 &\\
7 &1 &$u_{6}$ &$u_{6}$ &$\bar u_{6}$ &$u_{6}$ &$\bar u_{6}$ &1 &1 &$\bar u_{6}$ &$\bar u_{6}$ &1 &1 &$\bar u_{6}$ &$u_{6}$ &$\bar u_{6}$ &$u_{6}$ &$u_{6}$ &1 &0 &19 &\\
8 &1 &$u_{7}$ &$u_{1}$ &$\bar u_{4}$ &$u_{4}$ &$u_{8}$ &$u_{6}$ &$u_{3}$ &$u_{2}$ &$\bar u_{7}$ &$\bar u_{6}$ &$\bar u_{3}$ &$\bar u_{1}$ &$\bar u_{5}$ &$u_{5}$ &$\bar u_{8}$ &$\bar u_{2}$ &-1 &0 &19 &\\
9 &1 &$u_{8}$ &$\bar u_{4}$ &$\bar u_{2}$ &$u_{2}$ &$u_{4}$ &$\bar u_{6}$ &$u_{6}$ &$\bar u_{8}$ &$\bar u_{8}$ &$u_{6}$ &$\bar u_{6}$ &$u_{4}$ &$u_{2}$ &$\bar u_{2}$ &$\bar u_{4}$ &$u_{8}$ &1 &0 &19 &\\
10 &1 &-1 &-1 &1 &1 &1 &1 &-1 &1 &-1 &1 &-1 &-1 &-1 &-1 &1 &1 &-1 &0 &19 &\\
11 &1 &$\bar u_{8}$ &$u_{4}$ &$u_{2}$ &$\bar u_{2}$ &$\bar u_{4}$ &$u_{6}$ &$\bar u_{6}$ &$u_{8}$ &$u_{8}$ &$\bar u_{6}$ &$u_{6}$ &$\bar u_{4}$ &$\bar u_{2}$ &$u_{2}$ &$u_{4}$ &$\bar u_{8}$ &1 &0 &19 &\\
12 &1 &$\bar u_{7}$ &$\bar u_{1}$ &$u_{4}$ &$\bar u_{4}$ &$\bar u_{8}$ &$\bar u_{6}$ &$\bar u_{3}$ &$\bar u_{2}$ &$u_{7}$ &$u_{6}$ &$u_{3}$ &$u_{1}$ &$u_{5}$ &$\bar u_{5}$ &$u_{8}$ &$u_{2}$ &-1 &0 &19 &\\
13 &1 &$\bar u_{6}$ &$\bar u_{6}$ &$u_{6}$ &$\bar u_{6}$ &$u_{6}$ &1 &1 &$u_{6}$ &$u_{6}$ &1 &1 &$u_{6}$ &$\bar u_{6}$ &$u_{6}$ &$\bar u_{6}$ &$\bar u_{6}$ &1 &0 &19 &\\
14 &1 &$\bar u_{5}$ &$u_{7}$ &$u_{8}$ &$\bar u_{8}$ &$u_{2}$ &$u_{6}$ &$u_{3}$ &$\bar u_{4}$ &$u_{5}$ &$\bar u_{6}$ &$\bar u_{3}$ &$\bar u_{7}$ &$u_{1}$ &$\bar u_{1}$ &$\bar u_{2}$ &$u_{4}$ &-1 &0 &19 &\\
15 &1 &$\bar u_{4}$ &$u_{2}$ &$\bar u_{8}$ &$u_{8}$ &$\bar u_{2}$ &$\bar u_{6}$ &$u_{6}$ &$u_{4}$ &$u_{4}$ &$u_{6}$ &$\bar u_{6}$ &$\bar u_{2}$ &$u_{8}$ &$\bar u_{8}$ &$u_{2}$ &$\bar u_{4}$ &1 &0 &19 &\\
16 &1 &$\bar u_{3}$ &$\bar u_{3}$ &$\bar u_{6}$ &$u_{6}$ &$\bar u_{6}$ &1 &-1 &$\bar u_{6}$ &$u_{3}$ &1 &-1 &$u_{3}$ &$\bar u_{3}$ &$u_{3}$ &$u_{6}$ &$u_{6}$ &-1 &0 &19 &\\
17 &1 &$\bar u_{2}$ &$\bar u_{8}$ &$\bar u_{4}$ &$u_{4}$ &$u_{8}$ &$u_{6}$ &$\bar u_{6}$ &$u_{2}$ &$u_{2}$ &$\bar u_{6}$ &$u_{6}$ &$u_{8}$ &$u_{4}$ &$\bar u_{4}$ &$\bar u_{8}$ &$\bar u_{2}$ &1 &0 &19 &\\
18 &1 &$\bar u_{1}$ &$u_{5}$ &$\bar u_{2}$ &$u_{2}$ &$u_{4}$ &$\bar u_{6}$ &$\bar u_{3}$ &$\bar u_{8}$ &$u_{1}$ &$u_{6}$ &$u_{3}$ &$\bar u_{5}$ &$\bar u_{7}$ &$u_{7}$ &$\bar u_{4}$ &$u_{8}$ &-1 &0 &19 &\\

\hline
\end{tabular}
\label{tab.m19}
\end{table}

\begin{table}
\caption{$\chi_r(n)$, $m=20$, $\varphi(m)=8$.}
\begin{tabular}{r|cccccccccccccccccccc|c|r}
\hline
$r$ & 1 & 2 & 3 & 4 & 5 & 6 & 7 & 8 & 9 & 10 & 11 & 12 & 13 & 14 & 15 & 16 & 17 & 18 & 19 & 20 & $f$\\\hline
1 &1 &0 &1 &0 &0 &0 &1 &0 &1 &0 &1 &0 &1 &0 &0 &0 &1 &0 &1 &0 &1 &\\
2 &1 &0 &$-i$ &0 &0 &0 &$i$ &0 &-1 &0 &1 &0 &$-i$ &0 &0 &0 &$i$ &0 &-1 &0 &5 &\\
3 &1 &0 &-1 &0 &0 &0 &-1 &0 &1 &0 &1 &0 &-1 &0 &0 &0 &-1 &0 &1 &0 &5 &\\
4 &1 &0 &$i$ &0 &0 &0 &$-i$ &0 &-1 &0 &1 &0 &$i$ &0 &0 &0 &$-i$ &0 &-1 &0 &5 &\\
5 &1 &0 &-1 &0 &0 &0 &-1 &0 &1 &0 &-1 &0 &1 &0 &0 &0 &1 &0 &-1 &0 &4 &\\
6 &1 &0 &$i$ &0 &0 &0 &$-i$ &0 &-1 &0 &-1 &0 &$-i$ &0 &0 &0 &$i$ &0 &1 &0 &20 &\\
7 &1 &0 &1 &0 &0 &0 &1 &0 &1 &0 &-1 &0 &-1 &0 &0 &0 &-1 &0 &-1 &0 &20 &\\
8 &1 &0 &$-i$ &0 &0 &0 &$i$ &0 &-1 &0 &-1 &0 &$i$ &0 &0 &0 &$-i$ &0 &1 &0 &20 &\\

\hline
\end{tabular}
\label{tab.m20}
\end{table}

\small
\begin{table}
\caption{$\chi_r(n)$, $m=21$, $\varphi(m)=12$.}
\begin{tabular}{r|ccccccccccccccccccccc|c|r}
\hline
$r$ & 1 & 2 & 3 & 4 & 5 & 6 & 7 & 8 & 9 & 10 & 11 & 12 & 13 & 14 & 15 & 16 & 17 & 18 & 19 & 20 & 21 & $f$\\\hline
1 &1 &1 &0 &1 &1 &0 &0 &1 &0 &1 &1 &0 &1 &0 &0 &1 &1 &0 &1 &1 &0 &1 &\\
2 &1 &$u_{4}$ &0 &$\bar u_{4}$ &$\bar u_{2}$ &0 &0 &1 &0 &$u_{2}$ &$\bar u_{4}$ &0 &-1 &0 &0 &$u_{4}$ &$u_{2}$ &0 &$\bar u_{2}$ &-1 &0 &7 &\\
3 &1 &$\bar u_{4}$ &0 &$u_{4}$ &$\bar u_{4}$ &0 &0 &1 &0 &$u_{4}$ &$u_{4}$ &0 &1 &0 &0 &$\bar u_{4}$ &$u_{4}$ &0 &$\bar u_{4}$ &1 &0 &7 &\\
4 &1 &1 &0 &1 &-1 &0 &0 &1 &0 &-1 &1 &0 &-1 &0 &0 &1 &-1 &0 &-1 &-1 &0 &7 &\\
5 &1 &$u_{4}$ &0 &$\bar u_{4}$ &$u_{4}$ &0 &0 &1 &0 &$\bar u_{4}$ &$\bar u_{4}$ &0 &1 &0 &0 &$u_{4}$ &$\bar u_{4}$ &0 &$u_{4}$ &1 &0 &7 &\\
6 &1 &$\bar u_{4}$ &0 &$u_{4}$ &$u_{2}$ &0 &0 &1 &0 &$\bar u_{2}$ &$u_{4}$ &0 &-1 &0 &0 &$\bar u_{4}$ &$\bar u_{2}$ &0 &$u_{2}$ &-1 &0 &7 &\\
7 &1 &-1 &0 &1 &-1 &0 &0 &-1 &0 &1 &-1 &0 &1 &0 &0 &1 &-1 &0 &1 &-1 &0 &3 &\\
8 &1 &$\bar u_{2}$ &0 &$\bar u_{4}$ &$u_{4}$ &0 &0 &-1 &0 &$u_{2}$ &$u_{2}$ &0 &-1 &0 &0 &$u_{4}$ &$\bar u_{4}$ &0 &$\bar u_{2}$ &1 &0 &21 &\\
9 &1 &$u_{2}$ &0 &$u_{4}$ &$u_{2}$ &0 &0 &-1 &0 &$u_{4}$ &$\bar u_{2}$ &0 &1 &0 &0 &$\bar u_{4}$ &$\bar u_{2}$ &0 &$\bar u_{4}$ &-1 &0 &21 &\\
10 &1 &-1 &0 &1 &1 &0 &0 &-1 &0 &-1 &-1 &0 &-1 &0 &0 &1 &1 &0 &-1 &1 &0 &21 &\\
11 &1 &$\bar u_{2}$ &0 &$\bar u_{4}$ &$\bar u_{2}$ &0 &0 &-1 &0 &$\bar u_{4}$ &$u_{2}$ &0 &1 &0 &0 &$u_{4}$ &$u_{2}$ &0 &$u_{4}$ &-1 &0 &21 &\\
12 &1 &$u_{2}$ &0 &$u_{4}$ &$\bar u_{4}$ &0 &0 &-1 &0 &$\bar u_{2}$ &$\bar u_{2}$ &0 &-1 &0 &0 &$\bar u_{4}$ &$u_{4}$ &0 &$u_{2}$ &1 &0 &21 &\\

\hline
\end{tabular}
\label{tab.m21}
\end{table}

\normalsize
\small
\begin{table}
\caption{$\chi_r(n)$, $m=22$, $\varphi(m)=10$.}
\begin{tabular}{r|cccccccccccccccccccccc|c|r}
\hline
$r$ & 1 & 2 & 3 & 4 & 5 & 6 & 7 & 8 & 9 & 10 & 11 & 12 & 13 & 14 & 15 & 16 & 17 & 18 & 19 & 20 & 21 & 22 & $f$\\\hline
1 &1 &0 &1 &0 &1 &0 &1 &0 &1 &0 &0 &0 &1 &0 &1 &0 &1 &0 &1 &0 &1 &0 &1 &\\
2 &1 &0 &$\bar u_{2}$ &0 &$u_{4}$ &0 &$\bar u_{3}$ &0 &$\bar u_{4}$ &0 &0 &0 &$u_{1}$ &0 &$u_{2}$ &0 &$\bar u_{1}$ &0 &$u_{3}$ &0 &-1 &0 &11 &\\
3 &1 &0 &$\bar u_{4}$ &0 &$\bar u_{2}$ &0 &$u_{4}$ &0 &$u_{2}$ &0 &0 &0 &$u_{2}$ &0 &$u_{4}$ &0 &$\bar u_{2}$ &0 &$\bar u_{4}$ &0 &1 &0 &11 &\\
4 &1 &0 &$u_{4}$ &0 &$u_{2}$ &0 &$u_{1}$ &0 &$\bar u_{2}$ &0 &0 &0 &$u_{3}$ &0 &$\bar u_{4}$ &0 &$\bar u_{3}$ &0 &$\bar u_{1}$ &0 &-1 &0 &11 &\\
5 &1 &0 &$u_{2}$ &0 &$\bar u_{4}$ &0 &$\bar u_{2}$ &0 &$u_{4}$ &0 &0 &0 &$u_{4}$ &0 &$\bar u_{2}$ &0 &$\bar u_{4}$ &0 &$u_{2}$ &0 &1 &0 &11 &\\
6 &1 &0 &1 &0 &1 &0 &-1 &0 &1 &0 &0 &0 &-1 &0 &1 &0 &-1 &0 &-1 &0 &-1 &0 &11 &\\
7 &1 &0 &$\bar u_{2}$ &0 &$u_{4}$ &0 &$u_{2}$ &0 &$\bar u_{4}$ &0 &0 &0 &$\bar u_{4}$ &0 &$u_{2}$ &0 &$u_{4}$ &0 &$\bar u_{2}$ &0 &1 &0 &11 &\\
8 &1 &0 &$\bar u_{4}$ &0 &$\bar u_{2}$ &0 &$\bar u_{1}$ &0 &$u_{2}$ &0 &0 &0 &$\bar u_{3}$ &0 &$u_{4}$ &0 &$u_{3}$ &0 &$u_{1}$ &0 &-1 &0 &11 &\\
9 &1 &0 &$u_{4}$ &0 &$u_{2}$ &0 &$\bar u_{4}$ &0 &$\bar u_{2}$ &0 &0 &0 &$\bar u_{2}$ &0 &$\bar u_{4}$ &0 &$u_{2}$ &0 &$u_{4}$ &0 &1 &0 &11 &\\
10 &1 &0 &$u_{2}$ &0 &$\bar u_{4}$ &0 &$u_{3}$ &0 &$u_{4}$ &0 &0 &0 &$\bar u_{1}$ &0 &$\bar u_{2}$ &0 &$u_{1}$ &0 &$\bar u_{3}$ &0 &-1 &0 &11 &\\

\hline
\end{tabular}
\label{tab.m22}
\end{table}

Tables \ref{tab.m4}--\ref{tab.m22} are complemented to moduli $m\le 195$ in the ancillary file \texttt{anc/chi.txt}.
Each line of the file shows a single number $\chi_r(n)$ in the following format:
first $m$, then $r$, then $n$, then the absolute value $|\chi_r(n)|$ which is either zero or one,
a colon, and finally the argument of $\chi_r(n) = e^{i \pi \omega}$ represented by the fraction $\omega$.
Because $\chi_r(n)$ is completely multiplicative as a function of $n$ and also periodic in $n$,
the file enlists only terms of prime values $n$; the others may be obtained by prime factorization
of $n$.

\begin{rem}
Because the $\chi_r(n)$ are completely multiplicative, the sequences
of Dirichlet inverses are
\begin{equation}
\chi_r^{(-1)}(n) = \chi_r(n)\mu(n),
\label{eq.chiinv}
\end{equation}
where $\mu$ is the M\"obius function.
\end{rem}

\section{Dirichlet $L$-functions} \label{sec.L} 

\subsection{Symmetries}
The $L$-series are defined via the character table \cite{ShanksMC17,ApostolJNT2,BredeITSF17}.
\begin{defn} (Dirichlet $L$-functions)
\begin{equation}
L(s,\chi) = \sum_{n=1}^\infty \frac{\chi(n)}{n^s}
=\prod_{p}\frac{1}{1-\chi(p) p^{-s}}
=\frac{1}{m^s}\sum_{n=1}^m \chi(n)\zeta(s,n/m),
\label{eq.Lsdef}
\end{equation}
where $\zeta(.,.)$ is the Hurwitz zeta-function \cite{VepstasNA47,AdamchikAMC187} and the
product extends over all primes $p$.
\end{defn}
\begin{rem}
The Dirichlet inverse collects inverse powers of square-free $n$---rephrasing (\ref{eq.chiinv}):
\begin{equation}
\frac{1}{L(s,\chi)} = \prod_p [1-\chi(p)p^{-s}]
=
\sum_{n=1}^\infty \frac{\mu(n)\chi(n)}{n^s}
.
\end{equation}
\end{rem}

We provide an explicit table in Section \ref{sec.Ltab},
where the first column is the modulus $m$,
the second column the representation $r$ as
in chapter \ref{sec.char}, the third column $s$,
and the final
columns real and imaginary part of $L(s,\chi_r)$.

Cases where the values are complex conjugates of earlier values
have been replaced by fillers.  This happens whenever 
$\chi_r(n)=\bar \chi_{r'}(n)$ for all $n$.

Another symmetry
stems from \cite{Apostol}
\begin{equation}
L(s,\chi) = L(s,\psi) \prod_{p\mid m}\left(1-\frac{\psi(p)}{p^s}\right)
,
\end{equation}
where $d$ is an induced modulus, and $\psi(n)$ the character modulo $d$
identified by
\begin{equation}
\chi(n) = \psi(n)\chi_1(n)
.
\end{equation}
\begin{exa}
Examples are related  to conductors
$f<m$ in the tables of Section \ref{sec.char}:
\begin{itemize}
\item
$L(s,\chi_1)$ at $m=4$ equals $L(s,\chi_1)$ at $m=2$,
relating Tables \ref{tab.m4} and \ref{tab.m2}.
\item
$L(s,\chi_2)$ at $m=6$ equals $(1+2^{-s})L(s,\chi_2)$ at $m=3$,
relating Tables \ref{tab.m6} and \ref{tab.m3}.
\item
$L(s,\chi_1)$ at $m=8$ equals $L(s,\chi_1)$ at $m=2$
and
$L(s,\chi_3)$ at $m=8$ equals $L(s,\chi_2)$ at $m=4$,
relating Tables \ref{tab.m8}, \ref{tab.m4} and \ref{tab.m2}.
\item
$L(s,\chi_1)$ at $m=9$ equals $L(s,\chi_1)$ at $m=3$
and
$L(s,\chi_4)$ at $m=9$ equals $L(s,\chi_2)$ at $m=3$,
relating Tables \ref{tab.m9} and \ref{tab.m3}.
\item
$L(s,\chi_2)$ at $m=10$ equals $(1-i/2^s)L(s,\chi_2)$ at $m=5$,
relating Tables \ref{tab.m10} and \ref{tab.m5}.
\end{itemize}
\end{exa}

\subsection{Table}\label{sec.Ltab}
The list of $L(s,\chi_r)$ starts as follows:

\scriptsize
\normalsize

The case $s=1$ has been included for the non-principal characters,
where
\begin{equation}
\sum_{n=1}^m \chi_r(n)=0,\quad r>1
\label{eq.chisum}
\end{equation}
leads to cancellations which ensure
convergence \cite{HasseVorl,RaneProcIACS120}.
The
inner sum
in
\begin{equation}
L(s,\chi)= \sum_{k\ge 0}\sum_{n=1}^m \frac{\chi(n)}{(km+n)^s}
\end{equation}
is converted to a rational polynomial in $k$ with denominator degree
at least 2 larger than numerator degree, and the $k$-sum becomes essentially
an overlay
of Harmonic sums.
Table \ref{tab.chispec} shows the basic examples.
$\psi$ are the polygamma functions \cite[\S 6.4]{AS}.

\begin{table}
\caption{Closed forms of Dirichlet series.
}
\begin{tabular}{rrrcr}
\hline
$m$ & $r$  & $s$ & $L(s,\chi_r)$ & $L(s,\chi_r)$ \cite{EIS} \\
\hline
2 &  1 & $s$ & $(1-2^{-s})\zeta(s)$ & \\
2 &  1 & 2 & $\pi^2/8$ & A111003 \\
2 &  1 & 3 & $7\zeta(3)/8= -\psi''(1/4)/64-\pi^3/32$ & A233091 \\
2 &  1 & 4 & $\pi^4/96$ & \\
\hline
3 &  1 & $s$ & $(1-3^{-s})\zeta(s)$ & \\
3 &  1 & 2 & $4\pi^2/27$ & A214549 \\
3 &  1 & 3 & $26\zeta(3)/27=-\psi''(1/3)/27-4\pi^3/(81\surd 3)$ & \\
3 &  1 & 4 & $8\pi^4/729$ & A196751 \\
3 &  1 & 5 & $242\zeta(5)/243$ & \\
3 &  2 & 1 & $\pi/(3\surd 3)$ & A073010\\
3 &  2 & 2 & $2\psi'(1/3)/9-4\pi^2/27$ & A086724 \\
3 &  2 & 3 & $4\pi^3/(81\surd 3)$ & A129404 \\
3 &  2 & 4 & $\psi'''(1/3)/243-8\pi^4/729$ & \\
\hline
4 &  2 & 1 & $\pi/4$ & A003881 \\
4 &  2 & 2 & $\psi'(1/4)/8-\pi^2/8$ & A006752 \\
4 &  2 & 3 & $\pi^3/32$ & A153071 \\
4 &  2 & 4 & $\psi'''(1/4)/768-\pi^4/96$ & A175572 \\
\hline
5 &  1 & $s$ & $(1-5^{-s})\zeta(s)$ & \\
5 &  1 & 2 & $4\pi^2/25$ & \\
5 &  1 & 3 & $124\zeta(3)/125$ & \\
5 &  2 & 1 & $\frac{\pi}{5}[1/\sqrt{5-2\surd 5}+i/\sqrt{5+2\surd 5}]$ & \\
5 &  2 & 2 & $\frac{2}{25}[\psi'(1/5)-4\pi^2/(5-\surd 5)-i\psi'(2/5)+4i\pi^2/(5+\surd 5)]$ & \\
5 &  3 & 1 & $2\log[(1+\surd 5)/2]/\surd 5$ & A086466 \\
5 &  3 & 2 & $4 \pi^2/(25\surd 5)$ & \\
\hline
6 &  1 & $s$ & $(1+6^{-s}-2^{-s}-3^{-s})\zeta(s)$ & \\
6 &  1 & 2 & $\pi^2/9$ & A100044 \\
6 &  2 & 1 & $\pi/(2\surd 3)$ & A093766 \\
6 &  2 & 3 & $\pi^3/(18\surd 3)$ & \\
\hline
\end{tabular}
\label{tab.chispec}
\end{table}

\begin{rem}
For the principal character $\chi_1$, \cite{Apostol}
\begin{equation}
L(s,\chi_1)=\zeta(s)\prod_{p\mid m}\left(1-\frac{1}{p^s}\right),
\label{eq.Lchi1}
\end{equation}
so $L(s,\chi_1)$
is a rational multiple of $\pi^s$
if $s$ is even \cite{BradleyRJ6,TsumuraJNT48,HoffmanAMM102,AlkanMM163}.
\end{rem}

\subsection{First Derivative}
The derivative of (\ref{eq.Lsdef}) with respect to $s$ is
\begin{equation}
L'(s,\chi) = -\sum_{n\ge 2} \chi(n)\frac{\log n}{n^s}.
\end{equation}
For
$s>1$, values are calculated by direct application of (\ref{eq.Lsdef}),
\begin{equation}
L'(s,\chi)=\frac{1}{m^s}\left(-\log m
\sum_{n=1}^m\chi(n) \zeta(s,n/m)
+
\sum_{n=1}^m\chi(n) \zeta'(s,n/m)
\right).
\end{equation}
For $s=1$, each period of the character $\chi(n)=\chi(n+m)$ in
\begin{multline}
L'(1,\chi) = -\sum_{k\ge 0}\sum_{n=1}^m \chi(n)\frac{\log(km+n)}{km+n}
\\
=
-\sum_{n=1}^{\lfloor(m-1)/2\rfloor}
\chi(n)\frac{\log(n)}{n}
-\sum_{k=1}^\infty \sum_{n=-\lfloor(m-1)/2\rfloor}^{\lfloor(m-1)/2\rfloor} \chi(n)\frac{\log(km+n)}{km+n}
\end{multline}
is expanded in a series around a center
$X\equiv km$
to which
the indices have distances
$\epsilon=-n$,
\begin{multline}
\frac{\log(X-\epsilon)}{X-\epsilon} =
\frac{\log X}{X}+(-1+\log X)\frac{\epsilon}{X^2}+(-\frac{3}{2}+\log X)\frac{\epsilon^2}{X^3}
+(-\frac{11}{6}+\log X)\frac{\epsilon^3}{X^4}+\cdots
\\
=
\frac{\log X}{X}
+\sum_{j\ge 1}(-\sum_{u=1}^{j}\frac{1}{u}+\log X)\frac{\epsilon^j}{X^{j+1}}
.
\end{multline}
\begin{rem}
The generalization to other $s$ is \cite{MezoJNT130}
\begin{equation}
\frac{\log(X-\epsilon)}{(X-\epsilon)^s}
=
\frac{\log X}{X}+\sum_{j=1}^\infty \left[
-\sum_{u=0}^{j-1}\frac{(s)_u}{(j-u)u!}+\frac{(s)_j}{j!}\log X
\right] \frac{\epsilon^j}{X^{s+j}},
\end{equation}
where $(s)_j\equiv s(s+1)\cdots (s+j-1)$, $(s)_0=1$, is Pochhammer's symbol.
\end{rem}
Multiplication with $\chi(n)$ and summation over the classes $n$ for any non-principal character
cancels the leading $\log X/X$ term---see (\ref{eq.chisum})---and generates series \cite{LehmerAA27,IharaAA137}
\begin{equation}
L'(1,\chi) =
-\sum_{n=1}^{\lfloor(m-1)/2\rfloor}
\chi(n)\frac{\log(n)}{n}
-\sum_{k=1}^\infty \sum_{j\ge 2} \frac{\alpha_j+\beta_j\log(km)}{(km)^j}
\end{equation}
with constants $\alpha_j$ and $\beta_j$ depending on $\chi$ and on the order $j$.
\begin{rem}
Depending on the parity of the character, either the terms with even or those
with odd $j$ vanish.
This reduction of
terms has been driving the choice of $X$.
\end{rem}
Summations
\begin{equation}
\alpha_j\sum_{k \ge 1}\frac{1}{(km)^j}
=
\frac{\alpha_j}{m^j}\zeta(j)
,
\end{equation}
and
\begin{equation}
\beta_j\sum_{k \ge 1}\frac{\log(km)}{(km)^j}
=
\frac{\beta_j}{m^j}[\log m\,\zeta(j)-\zeta'(j)]
\end{equation}
rewrite $L'$ as a series over
Zeta Functions.
The $k$-terms up to some freely chosen $M\approx 40$
are summed directly to accelerate convergence:
\begin{alg}
\begin{multline}
L'(1,\chi)
=
-\sum_{n=2}^{Mm+\lfloor(m-1)/2\rfloor}
\chi(n)\frac{\log n}{n}
\\
-\sum_{j\ge 2}\left[
\frac{\alpha_j}{m^j}\zeta(j)+\frac{\beta_j}{m^j}[\zeta(j)\log m-\zeta'(j)]
-\sum_{1\le k\le M} \frac{\alpha_j+\beta_j \log(km)}{(km)^j}
\right]
.
\end{multline}
\end{alg}
\begin{rem}
This can probably be formulated 
much more efficiently in terms of the  Gauss sums \cite{IshibashiRM35,NicolCJM14}.
\end{rem}
The results for $L'(s,\chi_r)$ start as follows:

\scriptsize \begin{verbatim}
 m  r  s  Re(L')                                                Im(L')
 2  1  2 -0.41811583807616962590828560617781139711923775611996  0.00000000000000000000000000000000000000000000000000
 2  1  3 -0.06921016836272074748362627760768005927930796015366  0.00000000000000000000000000000000000000000000000000
 2  1  4 -0.01771623065878015673143625539935151371812832223282  0.00000000000000000000000000000000000000000000000000
 2  1  5 -0.00522011393348035186442563804164248718723344171968  0.00000000000000000000000000000000000000000000000000
 2  1  6 -0.00163309262853731373611967895541271800687857377431  0.00000000000000000000000000000000000000000000000000
 2  1  7 -0.00052595465151771844771643934267364011741770800196  0.00000000000000000000000000000000000000000000000000
 2  1  8 -0.00017197075267004922546732124391835435192533514509  0.00000000000000000000000000000000000000000000000000
 2  1  9 -0.00005669458151516296038034817251083265989289171083  0.00000000000000000000000000000000000000000000000000
 2  1 10 -0.00001877753891961746692900793448874128383181319749  0.00000000000000000000000000000000000000000000000000

 3  1  2 -0.63258236162649762708200575756500570930611425093721  0.00000000000000000000000000000000000000000000000000
 3  1  3 -0.14187732701989244492420427021033210559381975265229  0.00000000000000000000000000000000000000000000000000
 3  1  4 -0.05338083539288848748859404708635490112000898019216  0.00000000000000000000000000000000000000000000000000
 3  1  5 -0.02376820291781272030522239322060071116333138263190  0.00000000000000000000000000000000000000000000000000
 3  1  6 -0.01130138631860072648326657049664519466207000887400  0.00000000000000000000000000000000000000000000000000
 3  1  7 -0.00552422641475663255895936723479880503970450093810  0.00000000000000000000000000000000000000000000000000
 3  1  8 -0.00273338165371887933413456355075883756420758755113  0.00000000000000000000000000000000000000000000000000
 3  1  9 -0.00135998290214225322961491777490967865422189100126  0.00000000000000000000000000000000000000000000000000
 3  1 10 -0.00067839760411486935474878475316385767840903369166  0.00000000000000000000000000000000000000000000000000
 3  2  1  0.22266298696860150948666026276474436188657161605715  0.00000000000000000000000000000000000000000000000000
 3  2  2  0.13489092252345807409064587918196369202754729945570  0.00000000000000000000000000000000000000000000000000
 3  2  3  0.07527569062722762962331216028061202962234911131947  0.00000000000000000000000000000000000000000000000000
 3  2  4  0.04007380990190515691818068756656827444733394284007  0.00000000000000000000000000000000000000000000000000
 3  2  5  0.02075846279913862669323460796527515934928909079509  0.00000000000000000000000000000000000000000000000000
 3  2  6  0.01058516532249050924386643900827692965142388036618  0.00000000000000000000000000000000000000000000000000
 3  2  7  0.00534970119320939519232648443149160751998413890293  0.00000000000000000000000000000000000000000000000000
 3  2  8  0.00269034603907081869458536474978701888606095226299  0.00000000000000000000000000000000000000000000000000
 3  2  9  0.00134930467679590672574612091099997528531317044446  0.00000000000000000000000000000000000000000000000000
 3  2 10  0.00067573917582867185510738563249105312628698915310  0.00000000000000000000000000000000000000000000000000

 4  1  *  same block as m=2, r=1 above
 4  2  1  0.19290131679691242936318976402803278524509686762001  0.00000000000000000000000000000000000000000000000000
 4  2  2  0.08158073611659279510291216978594115145773887519886  0.00000000000000000000000000000000000000000000000000
 4  2  3  0.03157707945712738788724616744286091406753844890852  0.00000000000000000000000000000000000000000000000000
 4  2  4  0.01157054792451164165145272825991055112687726366004  0.00000000000000000000000000000000000000000000000000
 4  2  5  0.00409487498794859210069139346473920353166828166325  0.00000000000000000000000000000000000000000000000000
 4  2  6  0.00141739694744032189042296189559147079679885558589  0.00000000000000000000000000000000000000000000000000
 4  2  7  0.00048373355305791317333505432592425184709053568205  0.00000000000000000000000000000000000000000000000000
 4  2  8  0.00016362101836044233132872653651523994833718831247  0.00000000000000000000000000000000000000000000000000
 4  2  9  0.00005503463265438472183757623135734483380574563555  0.00000000000000000000000000000000000000000000000000
 4  2 10  0.00001844662524176733922410472153863013759986794730  0.00000000000000000000000000000000000000000000000000

 5  1  2 -0.79414955411761222469720192670858817205698401344150  0.00000000000000000000000000000000000000000000000000
 5  1  3 -0.18106414531976634265272520503820279224780330702467  0.00000000000000000000000000000000000000000000000000
 5  1  4 -0.06601391659734281254780710518639245085828845039340  0.00000000000000000000000000000000000000000000000000
 5  1  5 -0.02803059823039338648618168977463624264981206724792  0.00000000000000000000000000000000000000000000000000
 5  1  6 -0.01274655216161707783567191612822547752813257089256  0.00000000000000000000000000000000000000000000000000
 5  1  7 -0.00601266692474181544340142251056728527848680260418  0.00000000000000000000000000000000000000000000000000
 5  1  8 -0.00289780816429208681148594110523507826723163904303  0.00000000000000000000000000000000000000000000000000
 5  1  9 -0.00141515581506736113080511793801987035047486225809  0.00000000000000000000000000000000000000000000000000
 5  1 10 -0.00069686796644059205567056482789101956496211513721  0.00000000000000000000000000000000000000000000000000
 5  2  1  0.15455633174545896654216402988318223319631324458045 -0.04416511200957409045660097689544036344145251055763
 5  2  2  0.05050979313230396347463529389188839699537071768701 -0.06288371253648251957890661793245158835052207797768
 5  2  3  0.01507105948497285624292789584118242153531510296102 -0.04799063018512213874476530058592708250322082165337
 5  2  4  0.00424184607216138083790041032739339236339500612291 -0.03010137917480593897728507435161327269043871511918
 5  2  5  0.00114861565079184996476050087704459470274830844873 -0.01719585118829897996253763843920629701525900442239
 5  2  6  0.00030304940840236474636155422803615095905036896641 -0.00933236123888727278715072907924912021531877618753
 5  2  7  0.00007856512689696300450206984503991405022425728234 -0.00491427713632333387301179780383192804031806346576
 5  2  8  0.00002012751493845638650556527725292650875212564245 -0.00254037671848717154548440828698513402355480977756
 5  2  9  0.00000511524664696431121123658694448119840671360542 -0.00129802077872988575355809063127949978984266825936
 5  2 10  0.00000129298554161385699494858629450007762837834559 -0.00065830142225906966333365915041268268627304973615
 5  3  1  0.35624064703076149886468458637127319729600101706517  0.00000000000000000000000000000000000000000000000000
 5  3  2  0.20266211487080801527488733720785334085091900565061  0.00000000000000000000000000000000000000000000000000
 5  3  3  0.10388625869208176709713993389724045176780995165888  0.00000000000000000000000000000000000000000000000000
 5  3  4  0.05104273008990150830689821112259545031075078923703  0.00000000000000000000000000000000000000000000000000
 5  3  5  0.02473640090329101340316683838755544717970256852729  0.00000000000000000000000000000000000000000000000000
 5  3  6  0.01198055022616085816574664564068743667497734955752  0.00000000000000000000000000000000000000000000000000
 5  3  7  0.00582939152208435310011173652194507030486352146162  0.00000000000000000000000000000000000000000000000000
 5  3  8  0.00285323832859547683771660707120286577974119330454  0.00000000000000000000000000000000000000000000000000
 5  3  9  0.00140420989593554401404925170274317235732888543766  0.00000000000000000000000000000000000000000000000000
 5  3 10  0.00069416308492606798469708268390211537633278700735  0.00000000000000000000000000000000000000000000000000
 5  4  *  complex conjugate of block m=5, r=2 above

 6  1  2 -0.22106312441035038531715323839767664970922350123298  0.00000000000000000000000000000000000000000000000000
 6  1  3 -0.02384978528249851972724435786366761692582963218901  0.00000000000000000000000000000000000000000000000000
 6  1  4 -0.00373531726099356257429393164827326376542863699524  0.00000000000000000000000000000000000000000000000000
 6  1  5 -0.00065714165364763572457150480485955255715025474529  0.00000000000000000000000000000000000000000000000000
 6  1  6 -0.00012165894201189454533419008342071367659157573371  0.00000000000000000000000000000000000000000000000000
 6  1  7 -0.00002313970297458484939698379205383906267748210396  0.00000000000000000000000000000000000000000000000000
 6  1  8 -0.00000447269734236844559590073628663864560648820581  0.00000000000000000000000000000000000000000000000000
 6  1  9 -0.00000087354855781095376482688783796727343060881158  0.00000000000000000000000000000000000000000000000000
 6  1 10 -0.00000017180828814384318028884506434742083274351198  0.00000000000000000000000000000000000000000000000000
 6  2  1  0.12445616121617402941801461275592899266970070044133  0.00000000000000000000000000000000000000000000000000
 6  2  2  0.03322426198835216013333563169691854919063967469441  0.00000000000000000000000000000000000000000000000000
 6  2  3  0.00809032537282812029136942953937439659172001832609  0.00000000000000000000000000000000000000000000000000
 6  2  4  0.00185491362365415188068745267800728147981196264481  0.00000000000000000000000000000000000000000000000000
 6  2  5  0.00040825816321425536634349133557907194970873493671  0.00000000000000000000000000000000000000000000000000
 6  2  6  0.00008735165752376249519608704307916158530192375585  0.00000000000000000000000000000000000000000000000000
 6  2  7  0.00001832422685076895091380802807933370090496864408  0.00000000000000000000000000000000000000000000000000
 6  2  8  0.00000379090798028155362700114637172277474692570545  0.00000000000000000000000000000000000000000000000000
 6  2  9  0.00000077660160018611457020318183843244621675802787  0.00000000000000000000000000000000000000000000000000
 6  2 10  0.00000015799247627892726540465636999955610941416444  0.00000000000000000000000000000000000000000000000000

 7  1  2 -0.85309025126436133004946970313518822354904214298483  0.00000000000000000000000000000000000000000000000000
 7  1  3 -0.19072909719908906534161215116264945239873942439401  0.00000000000000000000000000000000000000000000000000
 7  1  4 -0.06800538708274220825519746102881664932138469374313  0.00000000000000000000000000000000000000000000000000
 7  1  5 -0.02845202516805281803177421465758424399556399041637  0.00000000000000000000000000000000000000000000000000
 7  1  6 -0.01283522907322728553125394548131939197516809981662  0.00000000000000000000000000000000000000000000000000
 7  1  7 -0.00603112705456884998744292659522213499402046341079  0.00000000000000000000000000000000000000000000000000
 7  1  8 -0.00290161312373374662114970629774909396918409617538  0.00000000000000000000000000000000000000000000000000
 7  1  9 -0.00141593387383728528835417591601881706412707398657  0.00000000000000000000000000000000000000000000000000
 7  1 10 -0.00069702611007132015872795136262345138110120062692  0.00000000000000000000000000000000000000000000000000
 7  2  1  0.14361034321910282672247043842346402006888045267441 -0.17369230259554111555029265006994697564320400789026
\end{verbatim}\normalsize

\section{Prime Zeta Modulo Functions} \label{sec.P} 

\subsection{Dirichlet Prime $L$-functions}
Restricting the summation of the Dirichlet series to primes defines
functions $S(s,\chi)$ that mediate between Dirichlet series and the Prime Zeta Function.
We recall an algorithmic standard of numerical computation
\cite{LanguascoMC78,Cohen,FlajoletVardi}.
The constituents are the numbers
\begin{equation}
l_p \equiv \frac{\chi(p)}{p^s}
\end{equation}
with Mercator series
\begin{equation}
-\log(1-\chi(p)/p^s)
=
-\log(1-l_p)
=\sum_{t\ge 1} \frac{1}{t} l_p^t.
\end{equation}
Because the Dirichlet inverse of $1/t$ is $\mu(t)/t$, the M\"obius inverse
is
\cite{LoxtonJAMS30,ErdelyiIII,WeisnerTAMS38}
\begin{equation}
l_p = -\sum_{t\ge 1} \frac{\mu(t)}{t}\log(1-[\chi(p)/p^s]^t)
\label{eq.lpmob}
.
\end{equation}
\begin{defn}
(Dirichlet Prime $L$-series)
\begin{equation}
S(s,\chi) \equiv \sum_p \frac{\chi(p)}{p^s}
.
\label{eq.Sdef}
\end{equation}
\end{defn}
To accelerate convergence, the sum is split at some threshold value $M$,
\begin{equation}
S(s,\chi) = \sum_{p\le M}\frac{\chi(p)}{p^s}+\sum_{p>M} \frac{\chi(p)}{p^s},
\quad
S(M,s,\chi) \equiv \sum_{p> M} \frac{\chi(p)}{p^s},
\end{equation}
in which (\ref{eq.lpmob}) is inserted at the right hand side:
\begin{equation}
S(s,\chi) = \sum_{p\le M}\frac{\chi(p)}{p^s}
-
\sum_{t\ge 1} \frac{\mu(t)}{t}
\sum_{p>M}
\log[1-\chi^t(p)p^{-st}]
\label{eq.Sstmp}
.
\end{equation}
Cutting also through the Euler products (\ref{eq.Lsdef})
\begin{equation}
L(s,\chi) = \prod_{p\le M}\frac{1}{1-\chi(p)p^{-s}}
\prod_{p> M}\frac{1}{1-\chi(p)p^{-s}}
\end{equation}
proposes
\begin{defn} (Incomplete Dirichlet Euler Product) \label{def.LMs}
\begin{equation}
L(M,s,\chi) \equiv L(s,\chi)\prod_{p\le M}(1-\frac{\chi(p)}{p^s}) = \prod_{p>M}\frac{1}{1-\chi(p)p^{-s}}
.
\end{equation}
\end{defn}
The logarithm
\begin{equation}
\log L(M,s,\chi) = -\sum_{p>M} \log(1-\chi(p)p^{-s})
\end{equation}
is inserted into the right hand side of (\ref{eq.Sstmp}).
\begin{alg}
\begin{equation}
S(s,\chi) = \sum_{p\le M}\frac{\chi(p)}{p^s}
+
\sum_{t\ge 1} \frac{\mu(t)}{t}
\log L(M,st,\chi^t).
\label{eq.algmu}
\end{equation}
\end{alg}
This generates a table of $S(s,\chi_r)$:

\scriptsize \begin{verbatim}
 m  r  s  Re(S)                                                 Im(S)
 2  1  2  0.20224742004106549850654336483224793417323134323989  0.00000000000000000000000000000000000000000000000000
 2  1  3  0.04976263929944353642311331466570670097541212192615  0.00000000000000000000000000000000000000000000000000
 2  1  4  0.01449313976424684494261929593315787016204105971484  0.00000000000000000000000000000000000000000000000000
 2  1  5  0.00450501748392425713281824253885571113169727672665  0.00000000000000000000000000000000000000000000000000
 2  1  6  0.00144508685063651295413367326605939920958594187454  0.00000000000000000000000000000000000000000000000000
 2  1  7  0.00047133285613359253512413872944872308918332888531  0.00000000000000000000000000000000000000000000000000
 2  1  8  0.00015515536651783056052343914268308052297714451207  0.00000000000000000000000000000000000000000000000000
 2  1  9  0.00005134257496245066307358514078311753682292034974  0.00000000000000000000000000000000000000000000000000

 3  1  2  0.34113630892995438739543225372113682306212023212878  0.00000000000000000000000000000000000000000000000000
 3  1  3  0.13772560226240649938607627762866966393837508488911  0.00000000000000000000000000000000000000000000000000
 3  1  4  0.06464746075190116593027361692081219114969538070250  0.00000000000000000000000000000000000000000000000000
 3  1  5  0.03163979114647569746203634953474048479424871705587  0.00000000000000000000000000000000000000000000000000
 3  1  6  0.01569834473815365973053970893135432376376975531762  0.00000000000000000000000000000000000000000000000000
 3  1  7  0.00782658548530597479392615061788036460724460003300  0.00000000000000000000000000000000000000000000000000
 3  1  8  0.00390898957624195798012410977216029436233090156130  0.00000000000000000000000000000000000000000000000000
 3  1  9  0.00195366231153715980294047535060885548327417269948  0.00000000000000000000000000000000000000000000000000
 3  2  1 -0.64194483853319570866139263972173431667541104401589  0.00000000000000000000000000000000000000000000000000
 3  2  2 -0.27470520828551848628957789816086063009994509882225  0.00000000000000000000000000000000000000000000000000
 3  2  3 -0.13052475556852058147112238237019273405338475461522  0.00000000000000000000000000000000000000000000000000
 3  2  4 -0.06372483064121439204992455625399328246981734399486  0.00000000000000000000000000000000000000000000000000
 3  2  5 -0.03151448029153218645003121115269230790129727571063  0.00000000000000000000000000000000000000000000000000
 3  2  6 -0.01568088471610726339651562337306385986464759578854  0.00000000000000000000000000000000000000000000000000
 3  2  7 -0.00782412274079501095457725272363147659937126191727  0.00000000000000000000000000000000000000000000000000
 3  2  8 -0.00390864007053523196268812155398435360012142057735  0.00000000000000000000000000000000000000000000000000
 3  2  9 -0.00195361255478026763866872767417697973059282808296  0.00000000000000000000000000000000000000000000000000

 4  1  * same block as m=2, r=1 above
 4  2  1 -0.33498132529999318106331712148754357377997538075508  0.00000000000000000000000000000000000000000000000000
 4  2  2 -0.09461989289295015794518679014917480960188034024972  0.00000000000000000000000000000000000000000000000000
 4  2  3 -0.03225247383350252743465978303821335087318799948530  0.00000000000000000000000000000000000000000000000000
 4  2  4 -0.01119397145618826174462577330588016650706523982717  0.00000000000000000000000000000000000000000000000000
 4  2  5 -0.00385806941548066209579442616677354317681180931505  0.00000000000000000000000000000000000000000000000000
 4  2  6 -0.00131658492330696537193003498997224601161685078059  0.00000000000000000000000000000000000000000000000000
 4  2  7 -0.00044569595893400199862208556511570721754211538618  0.00000000000000000000000000000000000000000000000000
 4  2  8 -0.00015003262315703762090694227805813265008262330846  0.00000000000000000000000000000000000000000000000000
 4  2  9 -0.00005031836933794511539693316106370484242184103002  0.00000000000000000000000000000000000000000000000000

 5  1  2  0.41224742004106549850654336483224793417323134323989  0.00000000000000000000000000000000000000000000000000
 5  1  3  0.16676263929944353642311331466570670097541212192615  0.00000000000000000000000000000000000000000000000000
 5  1  4  0.07539313976424684494261929593315787016204105971484  0.00000000000000000000000000000000000000000000000000
 5  1  5  0.03543501748392425713281824253885571113169727672665  0.00000000000000000000000000000000000000000000000000
 5  1  6  0.01700608685063651295413367326605939920958594187454  0.00000000000000000000000000000000000000000000000000
 5  1  7  0.00827103285613359253512413872944872308918332888531  0.00000000000000000000000000000000000000000000000000
 5  1  8  0.00405884536651783056052343914268308052297714451207  0.00000000000000000000000000000000000000000000000000
 5  1  9  0.00200395557496245066307358514078311753682292034974  0.00000000000000000000000000000000000000000000000000
 5  2  1  0.05558076757159122797438373110392029200397141552077  0.25744607231019845831553919953785168450933618650489
 5  2  2  0.00587745504741124399339243532481344811691373319858  0.15525350477709758586525784166633774678670632020307
 5  2  3  0.00061185910379617932279038303371260872339769646209  0.09055484241783400177149050749754610178341378435528
 5  2  4  0.00006064197279815635344055135957389107001014491352  0.05054452415058793994664192927792361208763038172014
 5  2  5  0.00000580003308992913574711743668771979341610600613  0.02719213801008999533003275544381628939733258482756
 5  2  6  0.00000054287198362827368021161276466778954793060719  0.01426158552668629463074677498250709051205907255207
 5  2  7  0.00000005018054931017382182274452449237991840398809  0.00735645310955466879930470745507895671903321075689
 5  2  8  0.00000000460549113371719109246140782054160828361855  0.00375400658116490237713537870580380742684581756294
 5  2  9  0.00000000042097057300122419349132192179194835304856  0.00190234443109101218028442898933440462363147712115
 5  3  1 -1.00799654793986117226166607551267856699903195664933  0.00000000000000000000000000000000000000000000000000
 5  3  2 -0.38071874461129711635917015923497145487567001845229  0.00000000000000000000000000000000000000000000000000
 5  3  3 -0.16474980378367446941888405281224241856194292947043  0.00000000000000000000000000000000000000000000000000
 5  3  4 -0.07523495147825830029432718746871438806203449993040  0.00000000000000000000000000000000000000000000000000
 5  3  5 -0.03542159890156654964644026157148756761652565948421  0.00000000000000000000000000000000000000000000000000
 5  3  6 -0.01700490923570389119218181500524906761727983414130  0.00000000000000000000000000000000000000000000000000
 5  3  7 -0.00827092778629868551875580553631885600410697969667  0.00000000000000000000000000000000000000000000000000
 5  3  8 -0.00405883591198651733558468768304964106526021571057  0.00000000000000000000000000000000000000000000000000
 5  3  9 -0.00200395472034918349581165226497195360585348682960  0.00000000000000000000000000000000000000000000000000

 6  1  2  0.09113630892995438739543225372113682306212023212878  0.00000000000000000000000000000000000000000000000000
 6  1  3  0.01272560226240649938607627762866966393837508488911  0.00000000000000000000000000000000000000000000000000
 6  1  4  0.00214746075190116593027361692081219114969538070250  0.00000000000000000000000000000000000000000000000000
 6  1  5  0.00038979114647569746203634953474048479424871705587  0.00000000000000000000000000000000000000000000000000
 6  1  6  0.00007334473815365973053970893135432376376975531762  0.00000000000000000000000000000000000000000000000000
 6  1  7  0.00001408548530597479392615061788036460724460003300  0.00000000000000000000000000000000000000000000000000
 6  1  8  0.00000273957624195798012410977216029436233090156130  0.00000000000000000000000000000000000000000000000000
 6  1  9  0.00000053731153715980294047535060885548327417269948  0.00000000000000000000000000000000000000000000000000
 6  2  1 -0.14194483853319570866139263972173431667541104401589  0.00000000000000000000000000000000000000000000000000
 6  2  2 -0.02470520828551848628957789816086063009994509882225  0.00000000000000000000000000000000000000000000000000
 6  2  3 -0.00552475556852058147112238237019273405338475461522  0.00000000000000000000000000000000000000000000000000
 6  2  4 -0.00122483064121439204992455625399328246981734399486  0.00000000000000000000000000000000000000000000000000
 6  2  5 -0.00026448029153218645003121115269230790129727571063  0.00000000000000000000000000000000000000000000000000
 6  2  6 -0.00005588471610726339651562337306385986464759578854  0.00000000000000000000000000000000000000000000000000
 6  2  7 -0.00001162274079501095457725272363147659937126191727  0.00000000000000000000000000000000000000000000000000
 6  2  8 -0.00000239007053523196268812155398435360012142057735  0.00000000000000000000000000000000000000000000000000
 6  2  9 -0.00000048755478026763866872767417697973059282808296  0.00000000000000000000000000000000000000000000000000

\end{verbatim}\normalsize

\begin{rem}
For the principal character $\chi_1$,
\begin{gather}
S(s,\chi_1) = P(s)-2^{-s}, \quad m=2,
\\
S(s,\chi_1) = P(s)-3^{-s}, \quad m=3,
\\
S(s,\chi_1) = P(s)-2^{-s}-3^{-s}, \quad m=6,
\\
S(s,\chi_1) = P(s)-\sum_{\substack{p\le m\\ (p,m)>1}}p^{-s},
\end{gather}
where $P(s)$ is the Prime Zeta Function \cite{FrobergBIT8,MerrifieldPRSL33,MatharArxiv0803}.
\end{rem}

\subsection{Prime Zeta Modulo Functions} \label{sec.pzmod}

The sum (\ref{eq.Sdef}) over all primes
may be divided into distinct residue classes
\begin{equation}
S(s,\chi) = \sum_{n=1}^m \sum_{p\equiv n \pmod m}\frac{\chi(p)}{p^s}
= \sum_{n=1}^m \chi(n) \sum_{p\equiv n \pmod m}\frac{1}{p^s}.
\label{eq.Ssplit}
\end{equation}
This defines non-overlapping
subseries
of the Prime Zeta Function.
\begin{defn} (Prime Zeta Modulo Functions)
\begin{equation}
P_{m,n}(s)\equiv 
\sum_{p\equiv n \pmod m} \frac{1}{p^s}
;\quad P_{m,n}(M,s)\equiv 
\sum_{\substack{p\equiv n \pmod m\\ p>M}} \frac{1}{p^s}
.
\label{eq.Pmndef}
\end{equation}
\end{defn}
\begin{rem}
The Prime Zeta-function $P(s)$ is recovered
if the filtering with
the modular combs is
undone:
\begin{equation}
\sum_{n=1}^m P_{m,n}(s) =  P(s).
\label{eq.psum}
\end{equation}
\end{rem}

The matrix of the characters is orthogonal
\cite{HasseVorl},
so the inversion of (\ref{eq.Ssplit}) is
\begin{equation}
P_{m,n}(s) = \frac{1}{\varphi(m)}\sum_{r=1}^{\varphi(m)} \bar \chi_r(n) S(s,\chi_r)
+\sum_{\substack{p=n\le m\\(p,m)>1}}\frac{1}{p^s},
\label{eq.inv}
\end{equation}
which acts on the individual terms and therefore remains valid for the incomplete sums,
\begin{equation}
P_{m,n}(M,s) = \frac{1}{\varphi(m)}\sum_{r=1}^{\varphi(m)} \bar \chi_r(n) S(M,s,\chi_r)
+\sum_{\substack{M<p=n\le m\\(p,m)>1}}\frac{1}{p^s}.
\end{equation}
The second terms take care of small primes for which $\chi(p)=0$.

The case $m=2$ is simply $P_{2,1}(s) = S(s,\chi_1)$ because there is only one even prime,
repeating values provided above.
The block for $m=4$, $s\le 8$ has been tabulated earlier \cite{SebahGourdon}.
Starting at $m\ge 3$, the $P_{m,n}(s)$ are:

\scriptsize \begin{verbatim}
 m  n  s  P
 3  1  2  0.03321555032221795055292717778013809648108756665327
 3  2  2  0.30792075860773643684250507594099872658103266547551
 3  1  3  0.00360042334694295895747694762923846494249516513694
 3  2  3  0.13412517891546354042859932999943119899587991975217
 3  1  4  0.00046131505534338694017453033340945433993901835382
 3  2  4  0.06418614569655777899009908658740273680975636234868
 3  1  5  0.00006265542747175550600256919102408844647572067262
 3  2  5  0.03157713571900394195603378034371639634777299638325
 3  1  6  0.00000873001102319816701204277914523194956107976454
 3  2  6  0.01568961472713046156352766615220909181420867555308
 3  1  7  0.00000123137225548191967444894712444400393666905787
 3  2  7  0.00782535411305049287425170167075592060330793097513
 3  1  8  0.00000017475285336300871799410908797038110474049198
 3  2  8  0.00390881482338859497140611566307232398122616106932
 3  1  9  0.00000002487837844608213587383821593787634067230826
 3  2  9  0.00195363743315871372080460151239291760693350039122
 3  1 10  0.00000000354755157005073612886511663730357082627733
 3  2 10  0.00097666493907697988040868523830634809488146371787

 4  1  2  0.05381376357405767028067828734153656228567550149509
 4  3  2  0.14843365646700782822586507749071137188755584174481
 4  1  3  0.00875508273297050449422676581374667505111206122043
 4  3  3  0.04100755656647303192888654885196002592430006070572
 4  1  4  0.00164958415402929159899676131363885182748790994383
 4  3  4  0.01284355561021755334362253461951901833455314977101
 4  1  5  0.00032347403422179751851190818604108397744273370580
 4  3  5  0.00418154344970245961430633435281462715425454302085
 4  1  6  0.00006425096366477379110181913804357659898454554698
 4  3  6  0.00138083588697173916303185412801582261060139632757
 4  1  7  0.00001281844859979526825102658216650793582060674956
 4  3  7  0.00045851440753379726687311214728221515336272213574
 4  1  8  0.00000256137168039646980824843231247393644726060181
 4  3  8  0.00015259399483743409071519071037060658652988391026
 4  1  9  0.00000051210281225277383832598985970634720053965986
 4  3  9  0.00005083047215019788923525915092341118962238068988
 4  1 10  0.00000010240775251510279580486929749957117091037360
 4  3 10  0.00001693866668446511506004583085023984513096217169

 5  1  2  0.01082089638114771753353951906172584388284719779619
 5  2  2  0.27586829355163944664905730184997372065557850052458
 5  3  2  0.12061478877454186078379946018363597386887218032151
 5  4  2  0.00494344133373647354014708373691239576593346459761
 5  1  3  0.00080913843084035641245250698022237496506614634497
 5  2  3  0.12815553197969650234624459561826033077604565502678
 5  3  3  0.03760068956186250057475408812071422899263187067150
 5  4  3  0.00019727932704417708966212394650976624166844988289
 5  1  4  0.00006986805789621433879330279589781606000671240287
 5  2  4  0.06292928488592025628255758548942987059983408077138
 5  3  4  0.01238476073533231633591565621150625851220369905124
 5  4  4  0.00000922608509805798535275143632392498999656748935
 5  1  5  0.00000625466213439143946805396018589577550095731367
 5  2  5  0.03131022310141769935983100374949396438572202646650
 5  3  5  0.00411808509132770402979824830567767498838944163893
 5  4  5  0.00000045462904446230372093652349817598208485130755
 5  1  6  0.00000056583972496957732807037158491679285049223691
 5  2  6  0.01563354178492824835195225955908066196274598028000
 5  3  6  0.00137195625824195372120548457657357145068690772793
 5  4  6  0.00000002296774134130364785875882024900330256162971
 5  1  7  0.00000005135773338184100299467054471296122828929121
 5  2  7  0.00781371671538540391312233979398137313283918252394
 5  3  7  0.00045726360583073511381763233890241641380597176705
 5  4  7  0.00000000117718407166718117192602022058130988530311
 5  1  8  0.00000000466637839516483023409561227013523337400965
 5  2  8  0.00390642361020853816259472105933508411048224883713
 5  3  8  0.00015241702904363578545934235353127668363643127419
 5  4  8  0.00000000006088726144763914163420444959362509039110
 5  1  9  0.00000000042413860329242757996461375187871653490432
 5  2  9  0.00195314978937341462986352384610597009748484035541
 5  3  9  0.00005080535828240244957909485677156547385336323426
 5  4  9  0.00000000000316803029120338647329183008676818185576
 5  1 10  0.00000000003855562522808473135592118445688884204375
 5  2 10  0.00097656604062943543015412120109166142062067593252
 5  3 10  0.00001693509508643674413086026760791168357555923310
 5  4 10  0.00000000000016548281548613787552698185521679533591

 6  1  2  0.03321555032221795055292717778013809648108756665327
 6  5  2  0.05792075860773643684250507594099872658103266547551
 6  1  3  0.00360042334694295895747694762923846494249516513694
 6  5  3  0.00912517891546354042859932999943119899587991975217
 6  1  4  0.00046131505534338694017453033340945433993901835382
 6  5  4  0.00168614569655777899009908658740273680975636234868
 6  1  5  0.00006265542747175550600256919102408844647572067262
 6  5  5  0.00032713571900394195603378034371639634777299638325
 6  1  6  0.00000873001102319816701204277914523194956107976454
 6  5  6  0.00006461472713046156352766615220909181420867555308
 6  1  7  0.00000123137225548191967444894712444400393666905787
 6  5  7  0.00001285411305049287425170167075592060330793097513
 6  1  8  0.00000017475285336300871799410908797038110474049198
 6  5  8  0.00000256482338859497140611566307232398122616106932
 6  1  9  0.00000002487837844608213587383821593787634067230826
 6  5  9  0.00000051243315871372080460151239291760693350039122
 6  1 10  0.00000000354755157005073612886511663730357082627733
 6  5 10  0.00000010243907697988040868523830634809488146371787

 7  1  2  0.00222617267552791635282763857623956021321204073965
 7  2  2  0.25309105534677419585838789294127430386526757184310
 7  3  2  0.11639715869790706431310433422397739230981168440572
 7  4  2  0.00918106219015021539397505080632887326456712985154
 7  5  2  0.04393429651968657897926715845093488956785137670029
 7  6  2  0.00700951134571340516000169799675822107497051929142
 7  1  3  0.00005804318475905574364933421290638511768635690102
 7  2  3  0.12510571610493702424985047895830983009774308084804
 7  3  3  0.03728334482215790002040591936793369431741171735180
 7  4  3  0.00076333454612616258831533213575487134206393917018
 7  5  3  0.00816313811942037718888138393185085116457285136558
 7  6  3  0.00047361062699928485358521008227468409628403051693
 7  1  4  0.00000175761045236349988868543636363067741669169080
 7  2  4  0.06250414531619014962877519174532221536444122236439
 7  3  4  0.01235886586608185233067861599145569967385803356856
 7  4  4  0.00006849185028472527342747805858054974145383323906
 7  5  4  0.00160798094502772328679680992735923678962260090859
 7  6  4  0.00003540504834664066899170677740848293815579997593
 7  1  5  0.00000005620336980628757219253533458064427574227765
 7  2  5  0.03125017021375850184161640531341302534198061658338
 7  3  5  0.00411596755817982752350156292375488984072031582867
 7  4  5  0.00000621245467572000204759207257104266311039915062
 7  5  5  0.00032040970697646805498369681388868538569365042382
 7  6  5  0.00000270232869773481537382030894090797347471275287
 7  1  6  0.00000000184792769612393206788740822662097980117679
 7  2  6  0.01562500714982673740450973299145410498133164686380
 7  3  6  0.00137178469979148069557289972220259075614886051984
 7  4  6  0.00000056453099043396973189225995347613645833557431
 7  5  6  0.00006402137112669812819433912919925890749360902649
 7  6  6  0.00000020739122115254537517376570565905253914018337
 7  1  7  0.00000000006176660482334123360593904477690007671484
 7  2  7  0.00781250030429406636325791106326665205045138947343
 7  3  7  0.00045724984468817721903313174844938465713849902148
 7  4  7  0.00000005131683551900888044022430562769825770021128
 7  5  7  0.00001280112105330664867218529101307058005909548744
 7  6  7  0.00000001594181701645953672715216978864714306104399
 7  1  8  0.00000000000208617300170230400726962584682245349808
 7  2  8  0.00390625001305502389361079137163542714377815396040
 7  3  8  0.00015241593481010996513706425970785646559051989722
 7  4  8  0.00000000466509238999135538389890295952403129187615
 7  5  8  0.00000256005892801890273373942320287012333528915764
 7  6  8  0.00000000122602055737564094051849217646524322084646
 7  1  9  0.00000000000007094349319404656244040521756222506490
 7  2  9  0.00195312500056290352684233687342252458639878703099
 7  3  9  0.00005080527189581143363371684961648541638328959879
 7  4  9  0.00000000042409795873240696030915167027969672594218
 7  5  9  0.00000051200309995434620611735051756683669451583231
 7  6  9  0.00000000009430265664074137638656701306377648898400
 7  1 10  0.00000000000000242353263829861631897758029628734812
 7  2 10  0.00097656250002434720384818986592574194488058629133
 7  3 10  0.00001693508830568579715275539844060030880701402880
 7  4 10  0.00000000003855433521437118429390189553064316830548
 7  5 10  0.00000010240016312422041208107787699079256289590033
 7  6 10  0.00000000000725388960799776561781675438249607954315

 8  1  2  0.00481719944001490430800016531982149659648111754763
 8  3  2  0.12380794753386495737887417239098394381326524063310
 8  5  2  0.04899656413404276597267812202171506568919438394745
 8  7  2  0.02462570893314287084699090509972742807429060111171
 8  1  3  0.00022482579077606892991207169983370296429249257040
 8  3  3  0.03795923737874621050097172756821732453455944389710
 8  5  3  0.00853025694219443556431469411391297208681956865002
 8  7  3  0.00304831918772682142791482128374270138974061680863
 8  1  4  0.00001240057914431765157022600362003113386760903516
 8  3  4  0.01242211735753787431164237341495396442453626263914
 8  5  4  0.00163718357488497394742653531001882069362030090868
 8  7  4  0.00042143825267967903198016120456505391001688713186
 8  1  5  0.00000071378969549987368070435289646431572140936406
 8  3  5  0.00412184874496031324554547883894169525654337148042
 8  5  5  0.00032276024452629764483120383314461966172132434174
 8  7  5  0.00005969470474214636876085551387293189771117154043
 8  1  6  0.00000004165020615431191777387933642274101526270831
 8  3  6  0.00137232803938853551236876963794880625582124542544
 8  5  6  0.00006420931345861947918404525870715385796928283867
 8  7  6  0.00000850784758320365066308449006701635478015090213
 8  1  7  0.00000000244227705387080249482274869134009349823858
 8  3  7  0.00045729980966012391873727903917860227855869182633
 8  5  7  0.00001281600632274139744853175941781659572710851098
 8  7  7  0.00000121459787367334813583310810361287480403030941
 8  1  8  0.00000000014348051242565858216560023784749195377568
 8  3  8  0.00015242051432548838564315112854960660265969307627
 8  5  8  0.00000256122819988404414966626671223608895530682613
 8  7  8  0.00000017348051194570507203958182099998387019083399
 8  1  9  0.00000000000843564128353088757035062026925167636757
 8  3  9  0.00005080569062402872344247706183871749797040231687
 8  5  9  0.00000051209437661149030743841950908607794886329229
 8  7  9  0.00000002478152616916579278208908469369165197837301
 8  1 10  0.00000000000049610802941223825264777808779980414291
 8  3 10  0.00001693512652591184376250291459388467539088986097
 8  5 10  0.00000010240725640707338356661664972148337110623069
 8  7 10  0.00000000354015855327129754291625635516974007231072

 9  1  2  0.00402158543754434812181742479406281668033541416633
 9  2  2  0.26039173226766023741732954512096620015351462748034
 9  4  2  0.00752255314281775931487450953309022672123401907520
 9  5  2  0.04312155344790480339285501400300902929359859108337
 9  7  2  0.02167141174185584311623524345298505307951813341174
 9  8  2  0.00440747289217139603232051681702349713391944691181
 9  1  3  0.00017008424576062139060622195078389739472568490955
 9  2  3  0.12580572311897707843586565043069790850018156355436
 9  4  3  0.00049409881828651507056510056350819662353118827975
 9  5  3  0.00810355558039841186646955334556152755964233043035
 9  7  3  0.00293624028289582249630562511494637092423829194765
 9  8  3  0.00021590021608805012626412622317176293605602576745
 9  1  4  0.00000825656380197190972707564864569710842392106615
 9  2  4  0.06256995660246944477532390284009835317875707673455
 9  4  4  0.00003616042894915105881679414653139348021542823421
 9  5  4  0.00160402376266128658885355483567653272045522630653
 9  7  4  0.00041689806259226397163066053823236375129966905345
 9  8  4  0.00001216533142704762592162891162785091054405930760
 9  1  5  0.00000041887811219089751742389446002587133878888430
 9  2  5  0.03125626271142129633931099332460795743618667830983
 9  4  5  0.00000272908573930689364034153551975616239291768163
 9  5  5  0.00032016550299262414302048445792287319178479768673
 9  7  5  0.00005950746362025771484480376104430641274401410669
 9  8  5  0.00000070750459002147370230256118556571980152038669
 9  1  6  0.00000002165314420330617700577923467863736900218271
 9  2  6  0.01562556625209709052016651772875887563938535907812
 9  4  6  0.00000020831510574072208249530862192022798415406849
 9  5  6  0.00006400699017674680116881995009262009505162631770
 9  7  6  0.00000850004277325413875254169128863308420792351334
 9  8  6  0.00000004148485662424219232847335759607977169015726
 9  1  7  0.00000000112936105133831008373964509933527361911788
 9  2  7  0.00781255137580477014545459089248034091018305787202
 9  4  7  0.00000001597315332810527813272756669841821856004569
 9  5  7  0.00001280029924401125735339208130937537206109959683
 9  7  7  0.00000121426974110247608623247991264625044448989429
 9  8  7  0.00000000243800171147144371869696620432106377350629
 9  1  8  0.00000000005916646720564647646539923300104797195776
 9  2  8  0.00390625466711536025997892590199064688088349496760
 9  4  8  0.00000000122706977496722943559128332841057231961576
 9  5  8  0.00000256001290170220771062441551706422585556920906
 9  7  8  0.00000017346661712083584208205240540896948444891846
 9  8  8  0.00000000014337153250371656534556461287448709689266
 9  1  9  0.00000000000310668372015285362010410932189277922708
 9  2  9  0.00195312542416744980289165249739994733511945265134
 9  4  9  0.00000000009433745499118703336858038475534132792061
 9  5  9  0.00000051200055837045246513736210734988438231726862
 9  7  9  0.00000002478093430737079598684953144379910656516057
 9  8  9  0.00000000000843289346544781165288562038743173047126
 9  1 10  0.00000000000016331196596546999409540120445937620322
 9  2 10  0.00097656253855672497852795390313933029665018147175
 9  4 10  0.00000000000725503565044170240033226878114294608674
 9  5 10  0.00000010240002421559270781200281133981669058688493
 9  7 10  0.00000000354013322243432895647068896731796850398736
 9  8 10  0.00000000000049603930917291933235567798154069536120

10  1  2  0.01082089638114771753353951906172584388284719779619
10  3  2  0.12061478877454186078379946018363597386887218032151
10  7  2  0.02586829355163944664905730184997372065557850052458
10  9  2  0.00494344133373647354014708373691239576593346459761
10  1  3  0.00080913843084035641245250698022237496506614634497
10  3  3  0.03760068956186250057475408812071422899263187067150
10  7  3  0.00315553197969650234624459561826033077604565502678
10  9  3  0.00019727932704417708966212394650976624166844988289
10  1  4  0.00006986805789621433879330279589781606000671240287
10  3  4  0.01238476073533231633591565621150625851220369905124
10  7  4  0.00042928488592025628255758548942987059983408077138
10  9  4  0.00000922608509805798535275143632392498999656748935
10  1  5  0.00000625466213439143946805396018589577550095731367
10  3  5  0.00411808509132770402979824830567767498838944163893
10  7  5  0.00006022310141769935983100374949396438572202646650
10  9  5  0.00000045462904446230372093652349817598208485130755
10  1  6  0.00000056583972496957732807037158491679285049223691
10  3  6  0.00137195625824195372120548457657357145068690772793
10  7  6  0.00000854178492824835195225955908066196274598028000
10  9  6  0.00000002296774134130364785875882024900330256162971
10  1  7  0.00000005135773338184100299467054471296122828929121
10  3  7  0.00045726360583073511381763233890241641380597176705
10  7  7  0.00000121671538540391312233979398137313283918252394
10  9  7  0.00000000117718407166718117192602022058130988530311
10  1  8  0.00000000466637839516483023409561227013523337400965
10  3  8  0.00015241702904363578545934235353127668363643127419
10  7  8  0.00000017361020853816259472105933508411048224883713
10  9  8  0.00000000006088726144763914163420444959362509039110
10  1  9  0.00000000042413860329242757996461375187871653490432
10  3  9  0.00005080535828240244957909485677156547385336323426
10  7  9  0.00000002478937341462986352384610597009748484035541
10  9  9  0.00000000000316803029120338647329183008676818185576
10  1 10  0.00000000003855562522808473135592118445688884204375
10  3 10  0.00001693509508643674413086026760791168357555923310
10  7 10  0.00000000354062943543015412120109166142062067593252
10  9 10  0.00000000000016548281548613787552698185521679533591

\end{verbatim}\normalsize

There are obvious refinements in the table whenever classes of
a (non-prime) modulus
are a union of
classes of smaller modulus
\cite{ZuckerPRSA464}:
\begin{exa}
The value $P_{6,1}(s)$ equals $P_{3,1}(s)$.
The value $P_{6,5}(s)$ equals $P_{3,2}(s)-1/2^s$.
\end{exa}
\begin{exa}
The sum $P_{8,1}(s)+P_{8,5}(s)$
equals $P_{4,1}(s)$.
The sum $P_{8,3}(s)+P_{8,7}(s)$
equals $P_{4,3}(s)$.
\end{exa}
\begin{exa}
The sum $P_{9,1}(s)+P_{9,4}(s)+P_{9,7}(s)$
equals $P_{3,1}(s)$.
\end{exa}
\begin{exa}
The $P_{m,n}(s)$ for the
four classes at $m=10$ are essentially a permutation of the four classes
at $m=5$, but $P_{5,2}(s)$ for $m=5$ is $P_{10,7}(s)+1/2^s$
at $m=10$.
\end{exa}

\subsection{Euler modulo product}
A further set of constants is obtained if the filtered
modular values of (\ref{eq.Pmndef}) are re-organized into associated Euler products.
\begin{defn} (Euler modulo product)
\begin{equation}
\zeta_{m,n}(s)\equiv \prod_{p=n \pmod m}\frac{1}{1-p^{-s}}
=
\sum_{q=1}^\infty \frac{c_{m,n}(q)}{q^s}.
\label{eq.zetadef}
\end{equation}
The completely multiplicative arithmetic function $c$ is defined as $c_{m,n}(q)=0$ if $q$ has at least one
prime factor $\neq n \pmod m$, and $c_{m,n}(q)=1$ if all prime factors
of $q$ are $\equiv n \pmod m$ or if $q=1$.
\end{defn}
\begin{rem}
The characteristic function $c_{m,n}(q)$ in (\ref{eq.zetadef})
selects integers $q$ that have prime power signatures with $k=1, 2,\ldots$ (not necessarily
distinct) primes
all in the same residuum class, which is isomorph to the combinatorics of almost-prime signatures
\cite{MatharArxiv0803}. Therefore
\begin{equation}
\zeta_{m,n}(s)=\sum_{k=0}^\infty P_{m,n,k}(s),
\end{equation}
where
\begin{equation}
P_{m,n,k}(s)\equiv \frac{1}{k!}\sum_{\substack{k_1+2k_2+\cdots kk_k=k\\ k_k\ge 0}}
(k;k_1k_2\cdots k_k)^* P_{m,n}^{k_1}(s)
P_{m,n}^{k_2}(2s)\cdots P_{m,n}^{k_k}(ks),
\end{equation}
are multinomials in the $P_{m,n}(s)$.
\end{rem}
The $\zeta_{m,n}$ are accessible via multiplicative mixing of the $L$-series
if
$\varphi(m)=2$
\cite{FlajoletVardi,FinchPAMS138}.
The
simpler approach implemented here is to accumulate
the
$P_{m,n}(s)$ obtained in Section \ref{sec.pzmod},
\begin{gather}
\log \zeta_{m,n}(s)
= -\sum_{p = n \pmod m}\log (1-p^{-s})
=
\sum_{t\ge 1}\frac{1}{t}P_{m,n}(st)
.
\end{gather}
This evaluation may again be split at $M$ using the same cut as in (\ref{eq.Pmndef}):
\begin{defn}
(Incomplete Euler modulo product)
\begin{equation}
\zeta_{m,n}(M,s)\equiv\prod_{\substack{p=n \pmod m\\p> M}}\frac{1}{1-p^{-s}}.
\end{equation}
\end{defn}
The associates to Definition \ref{def.LMs} are
\begin{gather}
\zeta_{m,n}(s)=
\zeta_{m,n}(M,s)
\prod_{\substack{p=n \pmod m\\p\le M}}\frac{1}{1-p^{-s}}
;
\\
\log\zeta_{m,n}(M,s)
=\sum_{t \ge 1}\frac{1}{t} P_{m,n}(M,st)
.
\label{eq.zetaM}
\end{gather}
\begin{rem}\label{rem.zetProd}
Reminiscent of (\ref{eq.psum}), the Riemann zeta function collects
all residue classes,
\begin{equation}
\prod_{n=1}^m \zeta_{m,n}(s) = \zeta(s).
\end{equation}
The product over the numbers listed times $\prod_{p\le m,(p,m)>1}(1-p^{-s})$
(that fall into the ``trenches'' where the $\chi=0$) compared with $\zeta(s)$
provides a check of numerical consistency.
\end{rem}

Since $\zeta_{2,2}(s)=1/(1-2^{-s})$ and $\zeta_{2,1}(s)=(1-2^{-s})\zeta(s)$,
the noteworthy values of $\zeta_{m,n}(s)$ start at $m=3$:

\scriptsize
\normalsize

Redundancies of $\zeta_{m,n}$ correlate with those of $P_{m,n}$ mentioned 
in Section \ref{sec.pzmod}---sums over $P$ replaced by products over $\zeta$:
\begin{exa}
$\zeta_{6,1}(s)$ equals $\zeta_{3,1}(s)$.
$\zeta_{6,5}(s)$ equals $(1-2^{-s})\zeta_{3,2}(s)$.
\end{exa}
\begin{exa}
$\zeta_{8,1}(s)\zeta_{8,5}(s)$ equals $\zeta_{4,1}(s)$.
$\zeta_{8,3}(s)\zeta_{8,7}(s)$ equals $\zeta_{4,3}(s)$.
\end{exa}
\begin{exa}
The $\zeta_{14,n}(s)$ are basically a permutation of the
$\zeta_{7,n}(s)$, but
$\zeta_{14,9}(s)=(1-2^{-s})\zeta_{7,2}(s)$.
\end{exa}
\begin{exa}
$\prod_{p=1 \pmod 4}(1+\frac{2}{p^3-1}) = \zeta_{4,1}^2(3)/\zeta_{4,1}(6)\approx 1.0176647145765659722043971208664301824136$
\cite[(1.12)]{AmdeberhanJNT128}.
\end{exa}

\section{Zeta Expansions} \label{sec.HL}
Euler products over the unrestricted set of primes split in that modulo basis 
via exponential product expansions
\cite{MoreeMM101}.
The exponents $\gamma_{s,j}^{(.)}$ and coefficients in the logarithmic series 
of the Euler products $A_1^{(s)}$, $Q_1^{(s)}$, $F_1^{(s)}$
and $C_1^{(s)}$ are exactly those of my earlier work
\cite{MatharArxiv0903}.
\begin{defn}
(Artin's constants of order $s$)
\begin{equation}
A_{m,n}^{(s)}\equiv \prod_{p=n \pmod m}\left(1-\frac{1}{p^s(p-1)}\right)
=
\prod_j \zeta_{m,n}(j)^{-\gamma_{s,j}^{(A)}}
.
\end{equation}
\end{defn}
\begin{equation}
\log A_{m,n}^{(s)} =
-\sum_{t=s+1}^\infty
P_{m,n}(t)
\sum_{j=1}^{\lfloor t/(1+s)\rfloor}
\frac{1}{j}\binom{t-sj-1}{j-1}
,
\end{equation}
In numerical practise we employ the variant
equivalent to (\ref{eq.zetaM})
with a threshold parameter $M$ 
defining $P_{m,n}(M,t)$
and defining an incomplete $A_{m,n}^{(s)}$.
The $A_{m,n}^{(s)}$ are:

\scriptsize \begin{verbatim}
 m  n  s  A
 3  1  1  0.96303628898369466607270650729810646818769403087295
 3  2  1  0.46597099348832591947626339239107570928665851405755
 3  *  1  0.37395581361920228805472805434641641511162924860615
 3  1  2  0.99586806888702059296181571677626946649772180127160
 3  2  2  0.74159506991523763012267987602493572177788195838723
 3  *  2  0.69750135849636590328467035082092292407315394621452
 3  1  3  0.99946588820812761768573268166616031108130297059616
 3  2  3  0.87316795462332634027719017321077110422790752841041
 3  *  3  0.85654044485354217442616798413595388216657280031765
 3  1  4  0.99992717965764321352395451288853559606744908398313
 3  2  4  0.93711768119261343014859982817460626921266125054119
 3  *  4  0.93126518416000433438923720555067698255842373458780
 3  1  5  0.99998983484977336049227878198437755796364841316094
 3  2  5  0.96867184697945078245775606279631864033868663395654
 3  *  5  0.96666886859677751274032837293001626421142382211819

 4  1  1  0.93618304689486840485374735217483602559545295571098
 4  3  1  0.79889464962976801798338935072077358715040862544174
 4  *  1  0.37395581361920228805472805434641641511162924860615
 4  1  2  0.98919985103315362579349439810395838096529088750735
 4  3  2  0.94015563220835776594586758756666336614933165088863
 4  *  2  0.69750135849636590328467035082092292407315394621452
 4  1  3  0.99794677819394517971323934178682919238458040490803
 4  3  3  0.98091740655598483179128987767545852498377442486971
 4  *  3  0.85654044485354217442616798413595388216657280031765
 4  1  4  0.99959625657107790051514518546892697831474365445595
 4  3  4  0.99375075010601300012201501946377971158220753835063
 4  *  4  0.93126518416000433438923720555067698255842373458780
 4  1  5  0.99991972912653683284175097874281623646737982713030
 4  3  5  0.99793184018617426954861166190925426257024950676203
 4  *  5  0.96666886859677751274032837293001626421142382211819

 5  1  1  0.98831956638305746883087282795162021735104254639947
 5  2  1  0.48531082947754999932252444786658812444646497522733
 5  3  1  0.82493234111986491690153714571045421637133702809175
 5  4  1  0.99485772777378307465830719513542008131033561396718
 5  *  1  0.37395581361920228805472805434641641511162924860615
 5  1  2  0.99911416667905889380339439800697914957598443799056
 5  2  2  0.74725940266248257162661615183905048245682166762209
 5  3  2  0.94387233155461696180996684651869263721075148724586
 5  4  2  0.99979302441165734257764415786782157555012208036948
 5  *  2  0.69750135849636590328467035082092292407315394621452
 5  1  3  0.99992325507141931782483437831240399829723649073316
 5  2  3  0.87456297038427688306528625223021965638998225414189
 5  3  3  0.98144009040883263601170435021757121554704060929700
 5  4  3  0.99999029508991648136798168771812469418297569009643
 5  *  3  0.85654044485354217442616798413595388216657280031765
 5  1  4  0.99999312300779961065284240241206938654053912931887
 5  2  4  0.93743420242628703757024971525968083790937519251505
 5  3  4  0.99382408909519320610177720172252856184425578311066
 5  4  4  0.99999952116182270606071177223396763347910848780771
 5  *  4  0.93126518416000433438923720555067698255842373458780
 5  1  5  0.99999937766961442465172997736994048856195463128477
 5  2  5  0.96874035025205510610892797235078265238823889538516
 5  3  5  0.99794215558580146608046689108068784979274640836406
 5  4  5  0.99999997579084470604995769973504106813807405835542
 5  *  5  0.96666886859677751274032837293001626421142382211819

 6  1  1  0.96303628898369466607270650729810646818769403087295
 6  5  1  0.93194198697665183895252678478215141857331702811511
 6  *  1  0.37395581361920228805472805434641641511162924860615
 6  1  2  0.99586806888702059296181571677626946649772180127160
 6  5  2  0.98879342655365017349690650136658096237050927784964
 6  *  2  0.69750135849636590328467035082092292407315394621452
 6  1  3  0.99946588820812761768573268166616031108130297059616
 6  5  3  0.99790623385523010317393162652659554768903717532618
 6  *  3  0.85654044485354217442616798413595388216657280031765
 6  1  4  0.99992717965764321352395451288853559606744908398313
 6  5  4  0.99959219327212099215850648338624668716017200057727
 6  *  4  0.93126518416000433438923720555067698255842373458780
 6  1  5  0.99998983484977336049227878198437755796364841316094
 6  5  5  0.99991932591427177544026432288652246744638620279385
 6  *  5  0.96666886859677751274032837293001626421142382211819

 7  1  1  0.99771564294568838736155996492787546367559477468412
 7  2  1  0.49840088692369542648778312810294075345529599481565
 7  3  1  0.82871803666235390545933514079733883693844062023419
 7  4  1  0.98998905141183604558273699292198729033368569094988
 7  5  1  0.94610323757016626868761072363193190084000156782121
 7  6  1  0.99248607820796097215463092567055639457362125429046
 7  *  1  0.37395581361920228805472805434641641511162924860615
 7  1  2  0.99994014189209871584452460006551005110995697243087
 7  2  2  0.74991747228874878297031470134449669622114285848704
 7  3  2  0.94419863275192592923695459672054797957947335493860
 7  4  2  0.99916135030829557805642100162223289908196026590518
 7  5  2  0.98983016697259050948663754611004452800101496357527
 7  6  2  0.99948806669231373616034944437904869460306577636950
 7  *  2  0.69750135849636590328467035082092292407315394621452
 7  1  3  0.99999818427485449426581627589239112109518799672518
 7  2  3  0.87499621737893938504603267043258995819277235099382
 7  3  3  0.98146776697270810153477930104022054022783560088610
 7  4  3  0.99992467473018342384264858111921255578780699414136
 7  5  3  0.99799160362258728065847744039609703225891846735173
 7  6  3  0.99996166797698051410762218911629473139618175755620
 7  *  3  0.85654044485354217442616798413595388216657280031765
 7  1  4  0.99999994188477665613746471458679347599374570540990
 7  2  4  0.93749983342357481295568615381996549141887642266209
 7  3  4  0.99382637891296178623822692780002794908822962790985
 7  4  4  0.99999316656592074861469998930147905723527287783743
 7  5  4  0.99959956791155045893890680467063608451168904075951
 7  6  4  0.99999707301010952296050938802699541446883588698900
 7  *  4  0.93126518416000433438923720555067698255842373458780
 7  1  5  0.99999999808814607332653001930682318795744533803823
 7  2  5  0.96874999276560345204224635032110372021185903647059
 7  3  5  0.99794234170960789630101905258256187373104983998247
 7  4  5  0.99999937902057396614860818298900423682749345076897
 7  5  5  0.99991997744742418905121707179384935036968868082300
 7  6  5  0.99999977533878029037234181445905570571069241006908
 7  *  5  0.96666886859677751274032837293001626421142382211819

\end{verbatim}\normalsize

The lines with a star in the $n$-column
are \cite{MatharArxiv0903}
\begin{equation}
A_1^{(s)}= \prod_{n=1}^m A_{m,n}^{(s)},
\end{equation}
i.e., the noted $A_{m,n}^{(s)}$
multiplied by the finite product
$\prod_p \{1-1/[p^s(p-1)]\}$
from primes up to and not coprime to $m$;
see Remark \ref{rem.zetProd}.
A005596 and A065414--A065416 in the OEIS display $A_1^{(s)}$ for $s\le 4$ \cite{EIS}.

\begin{defn}
(Quadratic Class numbers of order $s$)
\begin{equation}
Q_{m,n}^{(s)}\equiv \prod_{p=n \pmod m}\left(1-\frac{1}{p^s(p+1)}\right)
=
\prod_j \zeta_{m,n}(j)^{-\gamma_{s,j}^{(Q)}}
.
\end{equation}
\end{defn}
Again, the algorithm is the $M-$deferred variant
of a known formula \cite{MatharArxiv0903},
\begin{equation}
\log Q_{m,n}^{(s)} =
-\sum_{t=s+1}^\infty
P_{m,n}(t)
\sum_{j=1}^{\lfloor t/(1+s)\rfloor}
\frac{(-1)^{t-(s+1)j}}{j}\binom{t-sj-1}{j-1}
.
\end{equation}
The table of $Q_{m,n}^{(s)}$ starts:

\scriptsize \begin{verbatim}
 m  n  s  Q
 3  1  1  0.97024910467115917306803567015977644528885523914057
 3  2  1  0.79204649341104628208008707618094468295489090223018
 3  *  1  0.70444220099916559273660335032663721018858643141710
 3  1  2  0.99680763522216452381503254510333176771190098607201
 3  2  2  0.90960373636135175343904777594792044585740214764741
 3  *  2  0.88151383972517077692839182290322784712986925720808
 3  1  3  0.99959370433770807882786299647928332089368010624696
 3  2  3  0.95697939892931332309349292388556861626669448272441
 3  *  3  0.94773326214367537593952153765418961303363163231741
 3  1  4  0.99994499634458014625186893933369950019928710526941
 3  2  4  0.97889912348656360383240568800722815385319039104803
 3  *  4  0.97582415304766824167901143659479983197176497122921
 3  1  5  0.99999234838384640058266382093557284740952517853584
 3  2  5  0.98952999643109135827233602414981315609523742685183
 3  *  5  0.98850439774124690875110662385118666440095808327535

 4  1  1  0.95406039246755865186005688760123226145264727334098
 4  3  1  0.88603472890500783721446805186195985509972671129727
 4  *  1  0.70444220099916559273660335032663721018858643141710
 4  1  2  0.99262946073126248692585121240213394187141980334665
 4  3  2  0.96879198085650902105058006862070103808727094151243
 4  *  2  0.88151383972517077692839182290322784712986925720808
 4  1  3  0.99862038502725376088960590123790296746879346645464
 4  3  3  0.99030529616410002251835282806899969957047647793415
 4  *  3  0.94773326214367537593952153765418961303363163231741
 4  1  4  0.99973009395621554677063646518040509444080859765319
 4  3  4  0.99685542637431681147759855262917235702927491774292
 4  *  4  0.97582415304766824167901143659479983197176497122921
 4  1  5  0.99994643288800372098610092291896044680141586572085
 4  3  5  0.99896321876577456056788702369534663838265580294255
 4  *  5  0.98850439774124690875110662385118666440095808327535

 5  1  1  0.98994541651532837528981290049449971750949436487113
 5  2  1  0.81417661748357764967519511911124771377014221986475
 5  3  1  0.90845960510191820259198957462683801148959562318257
 5  4  1  0.99525215425415153950030018010344295622027317317902
 5  *  1  0.70444220099916559273660335032663721018858643141710
 5  1  2  0.99925503344071233662058074353202166787351604642839
 5  2  2  0.91411979597676293619735190404927935915076781045200
 5  3  2  0.97170968219087668252029279604183386856774563347235
 5  4  2  0.99981152134879912586741302748663655518106124260658
 5  *  2  0.88151383972517077692839182290322784712986925720808
 5  1  3  0.99993586794053861211626882911611188040468932785472
 5  2  3  0.95797248791315543969785734684555954155114677761234
 5  3  3  0.99070465607889613208223242762827729704972901822545
 5  4  3  0.99999120670666619470789435284342251314622703330412
 5  *  3  0.94773326214367537593952153765418961303363163231741
 5  1  4  0.99999426409769122005699815077922330151446689642523
 5  2  4  0.97911501946597037170638532422142794193185794775165
 5  3  4  0.99691092876378462606358844096686476050711455139439
 5  4  4  0.99999956721940822422487479489052376349550108528067
 5  *  4  0.97582415304766824167901143659479983197176497122921
 5  1  5  0.99999948124042658176935512678249171240320614070572
 5  2  5  0.98957593423045818490913570593415192280274622185561
 5  3  5  0.99897099455900934654360997881391321362270950221429
 5  4  5  0.99999997815156628432244203179623561003978532541884
 5  *  5  0.98850439774124690875110662385118666440095808327535

 6  1  1  0.97024910467115917306803567015977644528885523914057
 6  5  1  0.95045579209325553849610449141713361954586908267621
 6  *  1  0.70444220099916559273660335032663721018858643141710
 6  1  2  0.99680763522216452381503254510333176771190098607201
 6  5  2  0.99229498512147464011532484648864048638989325197899
 6  *  2  0.88151383972517077692839182290322784712986925720808
 6  1  3  0.99959370433770807882786299647928332089368010624696
 6  5  3  0.99858719888276172844538392057624551262611598197330
 6  *  3  0.94773326214367537593952153765418961303363163231741
 6  1  4  0.99994499634458014625186893933369950019928710526941
 6  5  4  0.99972676441180963795649942604993513585006678234693
 6  *  4  0.97582415304766824167901143659479983197176497122921
 6  1  5  0.99999234838384640058266382093557284740952517853584
 6  5  5  0.99994610165668179362257114019349539984360834713448
 6  *  5  0.98850439774124690875110662385118666440095808327535

 7  1  1  0.99783169567561415395612278331682981321793848654811
 7  2  1  0.83084436115736277012960608040205678524297640045329
 7  3  1  0.91204173657631900150492713248862784862731325200951
 7  4  1  0.99152664451130171888186444781116508688643378122073
 7  5  1  0.96301728477470669833755668029320880900891574734534
 7  6  1  0.99343749409772808721220289503432216846512348209610
 7  *  1  0.70444220099916559273660335032663721018858643141710
 7  1  2  0.99994366073240649656784328959429546613893530201614
 7  2  2  0.91657341209621074525122836181093861113002239080188
 7  3  2  0.97199490345740290514191139256429322419400292540344
 7  4  2  0.99929947053146742408555681415248114978845320558949
 7  5  2  0.99317882616912445691586599920122870731373256645422
 7  6  2  0.99955929237466022465905605429286339331226987661148
 7  *  2  0.88151383972517077692839182290322784712986925720808
 7  1  3  0.99999829680521359669791936140827355119253342304436
 7  2  3  0.95832951728972205581988628973822050125077925004215
 7  3  3  0.99072837048111798598099398660127807953178764545713
 7  4  3  0.99993720312563673896846596275692472208855362653198
 7  5  3  0.99865908524683796954174391736600314663890463047532
 7  6  3  0.99996710470497322840547153373281750802705003754425
 7  *  3  0.94773326214367537593952153765418961303363163231741
 7  1  4  0.99999994558480920546185156129258500018727263672918
 7  2  4  0.97916650671420237250529489743286478175232028721336
 7  3  4  0.99691288143967753375215937932656902663383764671792
 7  4  4  0.99999430503583323315981125802879383765208262150784
 7  5  4  0.99973294403625574083545067264002211232543320540325
 7  6  4  0.99999749025918214169988795654177516084234849072519
 7  *  4  0.97582415304766824167901143659479983197176497122921
 7  1  5  0.99999999821182133621213755516477429295696450571644
 7  2  5  0.98958332654672320978345933759407259691090563702888
 7  3  5  0.99897115320706590516862955106329442434349627356586
 7  4  5  0.99999948250950920843035994399632033679802298131632
 7  5  5  0.99994664636169323397746114726903151081859982442482
 7  6  5  0.99999980741214224538109413884142031379247371320710
 7  *  5  0.98850439774124690875110662385118666440095808327535

\end{verbatim}\normalsize
Equivalent to Remark \ref{rem.zetProd}, lines with a star in the $n$-column are
\begin{equation}
Q_1^{(s)} = \prod_{n=1}^m Q_{m,n}^{(s)}.
\end{equation}
The Niklasch values of
$Q_1^{(s)}$, $s\le 5$ are A065463 and A065465--A065468 in the OEIS \cite{EIS}.

\begin{defn}
(Feller-Tornier constants)
\begin{equation}
F_{m,n}^{(s)}\equiv \prod_{p=n \pmod m}\left(1-\frac{2}{p^s}\right)
=
\prod_j \zeta_{m,n}(j)^{-\gamma_{s,j}^{(F)}}
.
\end{equation}
\end{defn}
Based on \cite{MatharArxiv0903}
\begin{equation}
\log F_{m,n}^{(s)} =
-\sum_{t=1}^\infty
\frac{2^t}{t}
P_{m,n}(st)
,
\end{equation}
the $F_{m,n}^{(s)}$ become:

\scriptsize \begin{verbatim}
 m  n  s  F
 3  1  2  0.93484201367742708692711271624051049601103614106357
 3  2  2  0.44372767162332273561525453355445652782381952938988
 3  *  2  0.32263409893924467057953169254823706657095057966583
 3  1  3  0.99280761653561283079975439424575351687703956375940
 3  2  3  0.73634019702793989821414960486873306430411328482419
 3  *  3  0.67689273700988199361023732672438921279767839745979
 3  1  4  0.99907744600555800880829233535859351936265583252243
 3  2  4  0.87204973197830029518258016493629726761761662731638
 3  *  4  0.84973299138471876626505370362916043989282010424286
 3  1  5  0.99987468990135792609587644151500113697139680938750
 3  2  5  0.93688662911191901929922845845837418044555880822266
 3  *  5  0.92905919295966281511524587198420062376637612342100

 4  1  2  0.89484122456248817072566150690843732198754780892072
 4  3  2  0.72109797824075241583243117750350641933238009488227
 4  *  2  0.32263409893924467057953169254823706657095057966583
 4  1  3  0.98251462525135333541051144592244831503480272183851
 4  3  3  0.91858546035955424632246635972848794535685587324250
 4  *  3  0.67689273700988199361023732672438921279767839745979
 4  1  4  0.99670115119737736072418577538862522209341552226715
 4  3  4  0.97433761118740866145644826461858163613349229991035
 4  *  4  0.84973299138471876626505370362916043989282010424286
 4  1  5  0.99935305638694748001375333214108685197931481807727
 4  3  5  0.99163800636498496397333514320139495841461274016083
 4  *  5  0.92905919295966281511524587198420062376637612342100

 5  1  2  0.97845248073805828720735998931024731012225701761724
 5  2  2  0.47437087172424840070720711055620424815486422859867
 5  3  2  0.76307380442876918582110108389946175183811018177911
 5  4  2  0.99014350052786791479556620332141840911985643518766
 5  *  2  0.32263409893924467057953169254823706657095057966583
 5  1  3  0.99838190086287504193401769542494480352829930716181
 5  2  3  0.74526882528822037898411982038512181070115931122217
 5  3  3  0.92488231672061579498342311673595319971735637472649
 5  4  3  0.99960547324813363288137636241078869060175533783316
 5  *  3  0.67689273700988199361023732672438921279767839745979
 5  1  4  0.99986026431454149472891209002225661373275171174437
 5  2  4  0.87424877013139375867977083071682392372665769680115
 5  3  4  0.97523240905393793969475076372750984033219642558801
 5  4  4  0.99998154787827064113659669538698369167247356564844
 5  *  4  0.84973299138471876626505370362916043989282010424286
 5  1  5  0.99998749067686156347656227622782664876124713464261
 5  2  5  0.93738708184645276170019195939946207242296096965759
 5  3  5  0.99176387687679582568564123108609354980158058605822
 5  4  5  0.99999909074199348489740247879140513547360047968785
 5  *  5  0.92905919295966281511524587198420062376637612342100

 6  1  2  0.93484201367742708692711271624051049601103614106357
 6  5  2  0.88745534324664547123050906710891305564763905877975
 6  *  2  0.32263409893924467057953169254823706657095057966583
 6  1  3  0.99280761653561283079975439424575351687703956375940
 6  5  3  0.98178692937058653095219947315831075240548437976559
 6  *  3  0.67689273700988199361023732672438921279767839745979
 6  1  4  0.99907744600555800880829233535859351936265583252243
 6  5  4  0.99662826511805748020866304564148259156299043121872
 6  *  4  0.84973299138471876626505370362916043989282010424286
 6  1  5  0.99987468990135792609587644151500113697139680938750
 6  5  5  0.99934573771938028725251035568893245914192939543750
 6  *  5  0.92905919295966281511524587198420062376637612342100

 7  1  2  0.99555404713198948244184584044601450508178750229430
 7  2  2  0.49691435036587215810973241763330728415652207058013
 7  3  2  0.76957796537542257228306221725905383464002237343862
 7  4  2  0.98166945377655916042561037674324801490546565135125
 7  5  2  0.91277467867618400750325447073175642635692243929785
 7  6  2  0.98600841416269679563485056385653389004971689718240
 7  *  2  0.32263409893924467057953169254823706657095057966583
 7  1  3  0.99988391667262821518818274735430677265688389455317
 7  2  3  0.74984143188165689857204887699782457488954387149038
 7  3  3  0.92546983387929396966566578391620995434498396242975
 7  4  3  0.99847336720476224584831435899199360231273056600688
 7  5  3  0.98367895449893787275102115297201007664596764353105
 7  6  3  0.99905281257738254491306703845873477181422457136253
 7  *  3  0.67689273700988199361023732672438921279767839745979
 7  1  4  0.99999648478110131582478600860698357487906770714195
 7  2  4  0.87499274570389232686770562590711676511699533369021
 7  3  4  0.97528291952775784702284332764806085376404886702137
 7  4  4  0.99986301635151287548319500298252048598086571436272
 7  5  4  0.99678408919749714690368672145075838680245939386610
 7  6  4  0.99992918995830049683594411988434204666642097898894
 7  *  4  0.84973299138471876626505370362916043989282010424286
 7  1  5  0.99999988759326185799713244165714207306583554795699
 7  2  5  0.93749968084921148189646542628572899629520177899819
 7  3  5  0.99176807708490809250200263138026260692965359369641
 7  4  5  0.99998757509072907576280787952699016917459084613132
 7  5  5  0.99935918111048145887987957173392322340983918846234
 7  6  5  0.99999459534270191193438256531181007469985761332297
 7  *  5  0.92905919295966281511524587198420062376637612342100

\end{verbatim}\normalsize
The redundancies are
\begin{equation}
F_{6,1}^{(s)}=F_{3,1}^{(s)};
\quad
F_{6,5}^{(s)}= F_{3,2}^{(s)}/(1-2^{1-s}).
\end{equation}
Equivalent to Remark \ref{rem.zetProd}, lines with a star in the $n$-column are
\begin{equation}
F_1^{(s)} = \prod_{n=1}^m F_{m,n}^{(s)}
\end{equation}
as known \cite{MatharArxiv0903}.
$F_1^{(2)}$ is A065474 in the OEIS \cite{EIS}.

\begin{defn}
(Hardy-Littlewood constants)
\begin{equation}
C_{m,n}^{(s)}\equiv \prod_{\substack{p=n \pmod m\\ p>s}}\frac{p^{s-1}(p-s)}{(p-1)^s}
=
\prod_j \zeta_{m,n}(j)^{-\gamma_{s,j}^{(C)}}
.
\end{equation}
\end{defn}
Based on \cite{MatharArxiv0903}
\begin{equation}
\log C_{m,n}^{(s)} =
-\sum_{t=2}^\infty
\frac{s^t-s}{t}
P_{m,n}(t)
,
\end{equation}
the $C_{m,n}^{(s)}$ read:

\scriptsize \begin{verbatim}
 m  n  s  C
 3  1  2  0.95836482362808797487757575013966951613301317628078
 3  2  2  0.91845582471428180359847772409497399151465223444510
 3  *  2  0.66016181584686957392781211001455577843262336028473
 3  1  3  0.86751218171239491908907658476288886972026952686301
 3  2  3  0.73216995449044542201998526124057450315277708109614
 3  *  3  0.63516635460427120720669659127252241734206568733237
 3  1  4  0.72258658325534766847436578690250565003731801889725
 3  2  4  0.42554745117605448763271593111166763672980634588268
 3  *  4  0.30749487875832709312335448607107685302217851995066
 3  1  5  0.52346534507199690512570203587861939776333944516009
 3  2  5  0.78300290352910122201164489284179973132539677801906
 3  *  5  0.40987488508823647447878121233795527789635801325495
 3  1  6  0.27710605535004326975691888903900717559837928534397
 3  2  6  0.67343998355656740459077411005798448137406580764280
 3  *  6  0.18661429735835839665692484794418833784007394494559

 4  1  2  0.92306113221757975924356514047403417887101279477641
 4  3  2  0.71518753504535953339073535622822204345296583727742
 4  *  2  0.66016181584686957392781211001455577843262336028473
 4  1  3  0.74414958837024499127360452418064114273268267540782
 4  3  3  0.85354660478324400008551727753202067459718543554908
 4  *  3  0.63516635460427120720669659127252241734206568733237
 4  1  4  0.44091192969092295875998558949907810730433587749524
 4  3  4  0.69740657499079114022349189307772080129218829858526
 4  *  4  0.30749487875832709312335448607107685302217851995066
 4  1  5  0.83597682802260698833719482963834225780151289192106
 4  3  5  0.49029455284991749943661473521891278555391304309671
 4  *  5  0.40987488508823647447878121233795527789635801325495
 4  1  6  0.75229690708395736448771308785547959989749911962742
 4  3  6  0.24805937071004332290260100799138708454500603005911
 4  *  6  0.18661429735835839665692484794418833784007394494559

 5  1  2  0.98735247301453588313535598749142090390363595767527
 5  2  2  0.96641612578310894608587557568989440501197703436705
 5  3  2  0.74195287815639505497135252353507072030513398594048
 5  4  2  0.99464113053183188600569420316192917886306904257169
 5  *  2  0.66016181584686957392781211001455577843262336028473
 5  1  3  0.96012864968975466210115275000677088726813802997313
 5  2  3  0.89069345898428997391945247100168623128503858445870
 5  3  3  0.96661807398837421625616971838928560460691085028529
 5  4  3  0.98352374917967957309861313654634184087997254268265
 5  *  3  0.63516635460427120720669659127252241734206568733237
 5  1  4  0.91635769426465932563156516687748163064942656924781
 5  2  4  0.76395458708732551236931737716038770467396880305946
 5  3  4  0.93096047342619617143737099126866851277161150022652
 5  4  4  0.96628229818335611095281967822459833058807820685263
 5  *  4  0.30749487875832709312335448607107685302217851995066
 5  1  5  0.85407396102427097380601359543999477960804600701910
 5  2  5  0.57764567367339449979254761372677952826205357641920
 5  3  5  0.88138997545439239487651612116012354934561245726479
 5  4  5  0.94259730749543865248035195590178804888244553462061
 5  *  5  0.40987488508823647447878121233795527789635801325495
 5  1  6  0.77137630056008079159358507773745815248542947723585
 5  2  6  0.32451219145047009948254991056074596924905555805738
 5  3  6  0.81724710072367521785130485832746372286203725036689
 5  4  6  0.91220871378010997522334196835556787989433736577703
 5  *  6  0.18661429735835839665692484794418833784007394494559

 6  1  2  0.95836482362808797487757575013966951613301317628078
 6  5  2  0.91845582471428180359847772409497399151465223444510
 6  *  2  0.66016181584686957392781211001455577843262336028473
 6  1  3  0.86751218171239491908907658476288886972026952686301
 6  5  3  0.73216995449044542201998526124057450315277708109614
 6  *  3  0.63516635460427120720669659127252241734206568733237
 6  1  4  0.72258658325534766847436578690250565003731801889725
 6  5  4  0.42554745117605448763271593111166763672980634588268
 6  *  4  0.30749487875832709312335448607107685302217851995066
 6  1  5  0.52346534507199690512570203587861939776333944516009
 6  5  5  0.78300290352910122201164489284179973132539677801906
 6  *  5  0.40987488508823647447878121233795527789635801325495
 6  1  6  0.27710605535004326975691888903900717559837928534397
 6  5  6  0.67343998355656740459077411005798448137406580764280
 6  *  6  0.18661429735835839665692484794418833784007394494559

 7  1  2  0.99765398801824805745575199848586590104838977156875
 7  2  2  0.99668740154736185119549542357183826086432853703231
 7  3  2  0.74564040438722187807736297411080874332720005893822
 7  4  2  0.98906852604128092494920011866206688759944804504343
 7  5  2  0.93348561251953561042303182283359764377636008578062
 7  6  2  0.99193337190767667267067400622701760558938762002681
 7  *  2  0.66016181584686957392781211001455577843262336028473
 7  1  3  0.99284563576602975595901031078857636226019241830154
 7  2  3  0.98984306267452432665619052548686584276034684896281
 7  3  3  0.98203592018701854285981806936191667723888824054206
 7  4  3  0.96524528800661834895816280975502100995118818242669
 7  5  3  0.77093874395874628323816728886762758318909948223296
 7  6  3  0.97465778149976704833101678046677018357146007753252
 7  *  3  0.63516635460427120720669659127252241734206568733237
 7  1  4  0.98546514686054195689724549716779516136368771372367
 7  2  4  0.97925593342382577716321982515617176017498868620909
 7  3  4  0.96303959431570229715988416824794900218954678573364
 7  4  4  0.92635560049467742885049474387330246912217450891583
 7  5  4  0.47504942820632988734750925376986894672657306491884
 7  6  4  0.94695503795677660236507040844485582354536105333782
 7  *  4  0.30749487875832709312335448607107685302217851995066
 7  1  5  0.97541078884978790157546934143855452639378834459773
 7  2  5  0.96472640730111318180467286493968492162356472526847
 7  3  5  0.93672333482275622639932836356657031526367035173684
 7  4  5  0.86999884056567783769854473353741503126140009551379
 7  5  5  0.95367427149821135161097352295533167315911015618065
 7  6  5  0.90753546045655108160772300725983280806657638255495
 7  *  5  0.40987488508823647447878121233795527789635801325495
 7  1  6  0.96259053430958360955758267936940830333133617188597
 7  2  6  0.94607005665812543459253596932520179647066529165643
 7  3  6  0.90265831568425453492025154169848478262640590887244
 7  4  6  0.79353818470465233918261062748102769445314469036766
 7  5  6  0.92878975832084684493539660980925986457038262893686
 7  6  6  0.85504445986218770161405444452761140400159991710681
 7  *  6  0.18661429735835839665692484794418833784007394494559

\end{verbatim}\normalsize
For the standard reasons, $C_{6,1}^{(s)}=C_{3,1}^{(s)}$ and
$C_{6,5}^{(s)}=C_{3,2}^{(s)}$.
Lines with a star in the $n$-column are \cite{MatharArxiv0903}
\begin{equation}
C_1^{(s)} = \prod_{n=1}^m C_{m,n}^{(s)},
\end{equation}
compatible with constants A005597, A065418 and A065419 in the OEIS \cite{EIS}.

\bibliographystyle{amsplain}
\bibliography{all}

\end{document}